\RequirePackage{ifpdf}
\ifpdf % We are running pdfTeX in pdf mode
\documentclass[pdftex]{sigma}
\else
\documentclass{sigma}
\fi

\newtheorem{thm}{Theorem}[section]
\newtheorem{thm*}{Theorem}
\newtheorem{lem}[thm]{Lemma}
\newtheorem{prop}[thm]{Proposition}
\newtheorem{cor}[thm]{Corollary}

{\theoremstyle{definition}

\newtheorem{df}[thm]{Definition}

\newtheorem{rmk}[thm]{Remark}
\newtheorem{nota}[thm]{Notation}
}

\def\calB{{\mathcal B}}
\def\calP{{\mathcal P}}
\def\M{{\mathcal M}}
\def\Mco{{{}^{\rm c/o}\mathcal M}_{g,\delta_1,\dots,\delta_n}^{s,\beta}}
\def\Mcos{{{}^{\rm c/o}\mathcal M}_{g,\delta_1,\dots,\delta_n}^{s_1\dots s_k,\beta}}
\def\Sull{\rm {\mathcal S}ull^{\rm c/o}}
\def\wSull{\widetilde{\rm {\mathcal S}ull}^{\rm c/o}}
\def\Rp{{\mathbb R}_{>0}}
\def\la{\langle}
\def\ra{\rangle}

\def\idagger{r^{\dagger}}

\def\res{r}
\def\oto{\otimes \dots \otimes}
\newcommand{\aoa}[2]{a_{#1}\otimes \dots \otimes a_{#2}}

\def\wt{\omega}

\def\CO{{\mathcal O_\mu}}
\def\calO{{\mathcal O}}
\def\a{\alpha}
\def\t{\tau}
\def\B{\mathcal B}
\def\Z{\mathbb Z}
\def\R{\mathbb R}
\def\ba{\mathbf{a}}
\def\d{d^{{\rm CH}^*}}
\def\s{s^{{\rm CH}^*}}

\begin{document}

\allowdisplaybreaks

\renewcommand{\thefootnote}{$\star$}

\renewcommand{\PaperNumber}{036}

\FirstPageHeading

\ShortArticleName{Open/Closed Hochschild Actions}

\ArticleName{Open/Closed String Topology and Moduli Space\\ Actions
via Open/Closed Hochschild Actions\footnote{This paper is a contribution to the Proceedings of
the XVIIIth International Colloquium on Integrable Systems and Quantum
Symmetries (June 18--20, 2009, Prague, Czech Republic).  The full
collection is
available at
\href{http://www.emis.de/journals/SIGMA/ISQS2009.html}{http://www.emis.de/journals/SIGMA/ISQS2009.html}}}

\Author{Ralph M.~KAUFMANN}

\AuthorNameForHeading{R.M.~Kaufmann}

\Address{Department of Mathematics, Purdue University, 150 N.~University St.,\\
West Lafayette, IN 47907-2067, USA}
\Email{\href{mailto:rkaufman@math.purdue.edu}{rkaufman@math.purdue.edu}}
\URLaddress{\url{http://www.math.purdue.edu/~rkaufman/}}

\ArticleDates{Received October 30, 2009, in f\/inal form April 10, 2010;  Published online April 30, 2010}

\Abstract{In this paper we extend our correlation functions to the open/closed
case. This gives rise to actions of an open/closed version of the
Sullivan PROP as well as an action of the relevant moduli space.
There are several unexpected structures and conditions that arise in this extension
which are forced upon us by considering the open sector.
For string topology type operations, one cannot just consider graphs,
but has to take punctures into account and one has to restrict the
underlying Frobenius algebras. In the moduli space, one f\/irst has to
pass to a smaller moduli space which is closed under open/closed duality
and then consider  covers in order to account for the punctures.}


\Classification{55P48; 81T30; 57R30; 16E40; 55P50}

\renewcommand{\thefootnote}{\arabic{footnote}}
\setcounter{footnote}{0}

\section{Introduction}
There has been a lot of interest in studying open/closed theories from physics and mathematics.
In physics this goes back to boundary CFTs and D-branes with an extensive literature.  In mathematics
 motivation come from open/closed string topology, Gromov--Witten invariants and TQFT again
 with a virtual onslaught of ideas.  Sources relevant to our constructions are \cite{CS,C1,C2,sullsurv}.
In this spirit, there have been many interesting  forays into the subject of open/closed
operations \cite{TZ,godin,ramirez,HVZ,BCT}.

Our point of view comes from the geometry  provided in \cite{KP,hoch1} and the
operations on Hochschild complexes def\/ined in \cite{hoch2} via correlations functions.
In this paper, we extend these correlation functions  to the open/closed
case. This leads to a dg-action of  a Sullivan-type PROP yielding string topology type operations and a cell level moduli space action on Hochschild complexes.
There are several surprising details  and conditions which to our knowledge have not been fully
discussed previously. These obstacles make the passage from the closed case to the open/closed case
far from being evident.

\looseness=-1
 The f\/irst is that~-- unlike the closed case~-- in the open/closed case, the role of punctures cannot be
suppressed, as they can arise as the result of an open gluing. These punctures
 contribute new factors to the correlators.
Another consequence of the presence of punctures is that in the moduli space case the underlying ribbon
 graph needs
extra decorations marking possible punctures and one cannot def\/ine
the string topology type operations just by looking at open/closed
graphs of Sullivan type. One has to know the puncture structure
inside the complementary regions as well.

Secondly we need a compatibility equation for the dg-PROP to operate, which is
satisf\/ied if the coef\/f\/icient modules of the Hochschild complexes have a geometric origin.
A further unexpected detail is that the moduli space is not the f\/irst moduli space one would choose.
For each moduli space one has to pick out a subspace  that satisf\/ies open/closed duality and then consider
covers or rather spaces which are stratif\/ied by covers of the usual moduli spaces. These are brane labelled
 open/closed moduli spaces $\Mco$ of bordered surfaces with punctures on the boundary and  clusters of marked
 points on the interior.

The main results are
 \begin{thm*}
\label{theoremone}
 There is an  open/closed colored $\beta$ brane labelled c/o dg-PROP cell model of the open/closed
 colored brane labelled
topolgical quasi PROP $($see Appendix~{\rm \ref{quasiproppar})} $\wSull$ which acts
in a~brane labelled open/closed colored dg-PROP fashion on the
brane
 labelled Hochschild complexes for a $\beta$-Frobenius algebra which
 satisfies
the Euler condition $(E)$; see Definition~{\rm \ref{Eulerdef}}.
\end{thm*}

This theorem def\/ines open/closed string topology type actions for
compact manifolds that are simply connected.

\begin{thm*}
\label{theoremtwo}
There is  an operadic cell model associated to the $\beta$-brane labelled open/closed mo\-du\-li spaces $\Mco$
 which acts on $\beta$-labelled Hochschild co-chains via operadic correlation functions with values in
 a $\beta$-labelled $\rm Hom$ operad.
\end{thm*}

The actions in both cases are made possible by a discrete version of
the c/o action.
\begin{thm*}
\label{theoremthree} For a basic  $\B$-Frobenius algebra the c/o
structure of discretely weighted arc-graphs acts  on  he collection
of complexes $B(\beta)$ and the isomorphic Hochschild complexes via
the correlation functions $Y$ defined by equations \eqref{cordef}
and \eqref{cordefgraph}.
\end{thm*}

The restriction basic $\B$-Frobenius algebras is simply of expository nature. We can deal with general systems of $\calB$-Frobenius algebras. For this we would introduce a new
propagator formalism which we will do elsewhere as not to put even more technical
structures into this exposition.

We begin by reviewing the relevant structures from \cite{KP} in Section~\ref{reviewsection}.
We then go on to construct the relevant spaces which carry the topological structures.
The f\/irst of these is taken from \cite{KP} and is concerned with graphs on windowed surfaces. We def\/ine a generalization and a restriction of this structure.
The restriction is the space that yields moduli spaces of curves with marked points and tangent
vectors, while the extension is used to def\/ine the open/closed Sullivan PROP. We also brief\/ly
discuss the moduli space $\Mco$ which provides the chain models for the moduli space action.
In all these cases, we associate a chain complex to these spaces where each
basis element, or cell, is indexed by a graph of arcs on the given windowed surface.

In Section~\ref{geosection}, we review the open/closed gluing operations in the
geometrical setting. The algebraic counterpart is given by brane labelled bar
complexes of brane labelled systems of Frobenius algebras which we
introduce in Section~\ref{algsection}.

The correlators in the discrete case are given in Section~\ref{corsection}. This is the technical heart of the paper and
yields Theorem \ref{theoremthree}. Given a graph  on a windowed surface which has a discrete weighting on its edges or arc there is a  universal formula for a correlation function. This extends our formula of \cite{hoch2} to the open/closed case.
In Section~\ref{octopsection} we given the details of the PROPic version of the actions,
which lead to open/closed
string topology operations. The section culminates in Theorem~\ref{theoremone}.
 In Section~\ref{modulisection} we give the details on our moduli space actions and the proof of Theorem~\ref{theoremtwo}.
We close the main part of the paper with an outlook. In an appendix,
we recall the rather technical def\/inition of a brane labelled c/o
structure and the other operadic and PROPic structures we use.

\section{Review of the KP-model for open/closed strings}
\label{reviewsection}

We recall the main features of the KP-model for open/closed string
interactions via foliations~\cite{KP}.  The interactions are given by surfaces
with boundary and punctures together with a foliation. There are marked points
on the boundary (at least one per boundary) and possibly also marked points
in the surface. The part of the boundary between two marked points is called a window.

The foliation is thought of as being transverse to the propagating string and as keeping
track of splitting and recombining of pieces of string. The foliation and its partial transverse measure
are encoded in a graph of weighted arcs. Furthermore there is the data of a brane
labelling which keeps track of the branes that the strings might end on.

Given two windows, either on two disjoint surfaces or on the same
surface, we can glue them and the foliations together if their
weights agree on these windows. In \cite{KP} we showed that this
gluing gives rise to an c/o structure on the topological level and
induces  chain level operations, which descend to a bi-modular
operad on the homology level.

More technically, the setup is as follows:

\subsection{Arc graphs in brane labelled windowed surfaces}

A {\it windowed surface} $F=F^s_g(\delta_1,\ldots ,\delta _r)$ is a
smooth oriented surface of genus $g\geq 0$ with $s\geq 0$ punctures and
$r\geq 1$ boundary components together with the specif\/ication of a
non-empty f\/inite subset $\delta _i$ of each boundary component,
for $i=1,\ldots ,r$, and we let $\delta =\delta
_1\cup\cdots\cup\delta _r$ denote the set of all distinguished points in the boundary
$\partial F$ of $F$ and let  $\sigma$ denote the set of all punctures.
The set of components of
$\partial F-\delta$ is called the set $W$ of {\it windows}.

Furthermore one needs to specify a brane labelling $\beta$.
For this we f\/irst f\/ix a set $\calB$ of basic brane labels
and denote by $\calP(\calB)$  power set.  Notice that $\varnothing\in\calP(\calB) $,
this will encode the closed sector. The elements of $\calP(\calB)$ of cardinality
bigger than one should be thought of as intersecting branes. These
give room for extra data, but it is possible to set all the contributions
for these  ``higher intersection'' branes
to zero in a given model.

A {\it brane-labeling} on a windowed surface $F$ is a function
\[
\beta : \ \delta\coprod\sigma\to{\mathcal P}({\mathcal B}),
\]
where $\sqcup$ denotes the disjoint union, so that if $\beta
(p)=\varnothing$ for some $p\in\delta$, then $p$ is the unique point
of $\delta$ in its component of~$\partial F$.  A brane-labeling
may take the value $\varnothing$ at a puncture.

A window $w\in W$ on a windowed surface $F$ brane-labeled by
$\beta$ is called {\it closed} if the endpoints of $w$ coincide at
the point $p\in\delta$ and $\beta (p)=\varnothing$; otherwise, the
window $w$ is called {\it open}.
Each window def\/ines a pair of brane labels $\beta(w)$ which is
the pair $(S,T)$ of the brane labels of the beginning and end of the window (these
may coincide).

\subsubsection{Arc families}

We def\/ine the sets
\[
\delta
(\beta )=\{ p\in\delta :\beta (p)\neq\varnothing\},\qquad \sigma
(\beta )=\{ p\in\sigma :\beta (p)\neq\varnothing\} .
\]

Def\/ine a $\beta$-arc $a$ in $F$ to be an arc properly embedded in $F$ with its endpoints in $W$ so that~$a$ is not
homotopic, f\/ixing its endpoints, to $\partial F-\delta (\beta )$.  For example, given a distinguished point
$p\in\partial
F$, consider the arc lying in a small neighborhood that simply connects one side of $p$ to another in $F$; $a$ is
a $\beta$-arc if and only if
$\beta (p)\neq\varnothing$. We will call such an arc a {\it small arc} around $p$.

Two $\beta$-arcs are {\it parallel} if they are homotopic rel~$\delta$, and a $\beta$-arc family is the homotopy
class
rel~$\delta$ of a collection of $\beta$-arcs, no two of which are parallel.  Notice that we take homotopies rel~$\delta$
rather than
rel $\delta (\beta )$.

Given a positively weighted arc family in $F$, let us furthermore say that a window $w\in W$ is {\it active}
if there
is an arc in the family with an endpoint in $w$, and otherwise the window is {\it inactive}.

\subsubsection{The mapping class group and arc graphs}

The {\it $($pure$)$ mapping class group} $MC(F)$ of $F$ is the group of
orientation-preserving homeomorphisms of $F$ pointwise f\/ixing
$\delta\cup\sigma$ modulo homotopies pointwise f\/ixing $\delta\cup\sigma$.

$MC(F)$ acts naturally on the set of $\beta$-arc families.
An {\it arc graph} is an equivalence graph under this action.

\subsection{Arc spaces, moduli spaces and the open/closed Sullivan PROP}
\subsubsection{The Arc spaces}

A {\it  weighting} on an arc family is the assignment of
a positive  real number to each of its components.
A weighting naturally passes to the arc graph.
A weighting is called {\it discrete} if it takes values in ${\mathbb N}$.

Let ${\rm Arc}'(F,\beta )$ denote the geometric realization of the partially
ordered set of all $\beta$-arc families in $F$.  ${\rm Arc}'(F,\beta )$ is
described as the set of all projective positively weighted $\beta$-arc
families in $F$ with the natural topology. (See for instance \cite{KLP} or \cite{P1} for
further detail.)

$MC(F)$ again
acts naturally on ${\rm Arc}'(F,\beta )$.  The {\it arc
complex} is def\/ined to be the  quotient under this action
\[
{\rm Arc}(F,\beta )={\rm Arc}'(F,\beta )/MC(F).
\]

We shall also consider the corresponding deprojectivized versions:
$\widetilde{\rm Arc}{}'(F,\beta )\approx {\rm Arc}'(F,\beta )\times{\mathbb R}_{>0}$ is the space of all positively weighted arc
families in $F$ with the natural topology, and
\[
\widetilde{\rm Arc}(F,\beta )=\widetilde{\rm Arc}{}'(F,\beta )/MC(F)\approx
{\rm Arc}(F,\beta )\times{\mathbb R}_{>0}.
\]

For any windowed surface $F$, def\/ine
\[
\widetilde{\rm Arc}(F)=\bigsqcup \widetilde{\rm Arc}(F,\beta),
\]
where the disjoint union is over all brane-labellings on $F$.
\[
\widetilde{\rm Arc}(n,m)=\bigsqcup \left\{ \alpha\in \widetilde{\rm Arc}(F):
\begin{array}{l}\alpha ~{\rm has}~n~{\rm closed~and}~ m~{\rm
open~active}\\{\rm
windows~and~no~inactive~windows}\end{array}\right\},
\] where the
disjoint union is over all orientation-preserving
homeomorphism classes of windowed surfaces.

\subsection{The open/closed Sullivan spaces}
As proved in \cite{KP}, the spaces $\widetilde{\rm Arc}(n,m)$ form an
c/o structure, see the appendix for the complete def\/inition. We will
now construct a suitable PROP to capture string topology type
operations. First, we have to add additional data to the surface
$(F,\beta)$, which is a partitioning ${\rm i/o} :=\{W_{\rm in},W_{\rm out}\}$ of
all of the windows of $F$ into ``in'' and  ``out'' windows. A
$\beta$-arc  family (or arc graph) on such a surface is called of
{\em Sullivan type} if
\begin{enumerate}\itemsep=0pt
\item arcs only run from ``in'' windows to ``out'' windows,
\item all {\em in} windows are active.
\end{enumerate}

We set
${\Sull}'(F,\beta,{\rm i/o})$ the geometric realization of the partially ordered set of all $\beta$-arc families
of Sullivan type. And let  ${\Sull}(F,\beta,{\rm i/o})$ be the quotient under the action of $MC(F)$.
As above, we also consider the deprojectivized versions $\wSull(F,\beta,{\rm i/o})$
\[
\wSull(F)=\bigsqcup \wSull(F,\beta,{\rm i/o}),
\]
where the disjoint union is over all brane-labellings on $F$ and partitions ${\rm i/o}$.
\begin{gather*}
\wSull(n_1,n_2,m_1,m_2)= \bigsqcup \Bigg\{ \alpha\in\Sull(F):  \\
\hphantom{\wSull(n_1,n_2,m_1,m_2)= \bigsqcup \Bigg\{}{}
\left. \begin{array}{@{}l@{}}\alpha ~{\rm has}~n=n_1+n_2~{\rm closed~and}~ m=m_1+m_2~{\rm open}\\
\text{windows with $n_1$ and $m_1$  active closed resp.\ open ``in''}\\
\text{windows and $n_2$, $m_2$  closed resp.\ open ``out'' windows}
\end{array}\right\},
\end{gather*}
where the
disjoint union is over all orientation-preserving
homeomorphism classes of windowed surfaces.

Notice that $\wSull(n,m)$ is also graded  by the number $k$ of
inactive ``out'' windows~-- by def\/inition all ``in'' windows are active.

\subsection{Moduli space} The spaces of $\beta$ arc families contains a moduli space.
We call an arc family or graph {\it quasi-filling}, if all complementary regions
are polygons or once punctured polygons.

The weighted quasi-f\/illing graphs form a subspace of each ${\rm Arc}(F,\beta)$ which
is homeomorphic to the moduli spaces $\M_{g,t_1,\dots, t_n,}^s$ of genus $g$ surfaces with $n+s$
marked points, where for the f\/irst $n$ marked points $t_i\geq1$ tangent vectors at the $i$th marked
point are specif\/ied. This is independent of $\beta$ since the small arcs only thicken the moduli
spaces by a factor of $\Rp$.
This fact is straightforward using the dual graph and Strebel dif\/ferentials~\cite{St}.

\subsubsection[Open/closed duality moduli space of open/closed  brane
labelled surfaces with marked point clusters]{Open/closed duality moduli space of open/closed  brane
labelled surfaces\\ with marked point clusters}

The moduli space  $\M_{g,t_1,\dots, t_n,}^s$  is actually too small and too big
for its cell model to act, as we will discuss below. It is too big in the sense that there are gluings
which when allowing open windows on boundaries unexpectedly take us out of moduli space on the chain level; see Figs.~\ref{example} and~\ref{example2}.
 There is however  a subspace of this space  given by the arc families
that are in general position with respect to the open/closed duality (see Fig.~\ref{openclosed}) which solves this problem.
The other problem that arises which is unique to the open sector is that unlike in~\cite{hoch1,hoch2}
we cannot restrict to the case of no punctures; $s=0$. But gluing with internal punctures is again not stable in the moduli space case. This problem is overcome by introducing clusters of brane labelled points. The resulting space is the open/closed duality moduli space of open/closed brane
labelled  surfaces with marked point clusters $\Mco$.
The details are given in Section~\ref{modulisection} below.

\section{Geometric c/o structures}
\label{geosection}

In \cite{KP} we axiomatized a structure of spaces that allow four dif\/ferent types of gluing: i.e.\
either open or closed and either self or non-self gluings together with a coloring of the open gluings by brane labels~$\beta$. In the above case the gluing is on two  windows, which are either both
open or closed and either belong to the same (self gluing) or disjoint surfaces (non-self gluing).
The coloring is by the brane labels. The color of a window is the pair of brane labels of its
end points. When we glue $\beta$ arc families, we glue the surfaces and
glue the arc families as foliations. We will now review this process according to~\cite{KP}.

\subsection{The gluing underlying the topological c/o structure}

If $\alpha\in\widetilde{\rm Arc} (n,m)$, then def\/ine the {\it $\alpha$-weight} $\alpha(w)$ of an active window $w$ to be the sum
of
the weights of arcs in $\alpha$ with endpoints in $w$, where we count with multiplicity
(so if an arc in $\alpha$ has both endpoints in $w$, then the weight of this arc contributes twice to the weight of
$w$).

Suppose we have a pair of arc families $\alpha_1$, $\alpha _2$
in respective windowed surfaces $F_1$, $F_2$ and a~pair of active windows $w_1$ in $F_1$ and $w_2$ in $F_2$, so that the
$\alpha _1$-weight of $w_1$ agrees with the $\alpha _2$-weight of $w_2$.
Since $F_1$, $F_2$ are oriented surfaces, so are the windows $w_1$, $w_2$ oriented.   In each operation, we
identify windows reversing orientation, and we identify certain distinguished points.

To def\/ine the open and closed gluing ($F_1\neq F_2$) and self-gluing ($F_1=F_2$) of $\alpha_1$, $\alpha _2$ along the
windows
$w_1$, $w_2$, we
identify windows and distinguished points in the natural way and combine foliations.

A crucial dif\/ference  between the closed and open string operations is that in the closed case,
the points are thought of as marked, which in the open case the points behave like punctures.
This means that in the closed case, we
replace the  distinguished point by simply forgetting that is was distinguished. This way  no puncture is created. In the open case the distinguished
points
always give rise to other distinguished points or perhaps punctures.
 In any case whenever distinguished points are
identif\/ied, one takes the union of brane labels (the intersection of branes) at the new resulting distinguished point or puncture.

\begin{figure}[th]
\centerline{\includegraphics{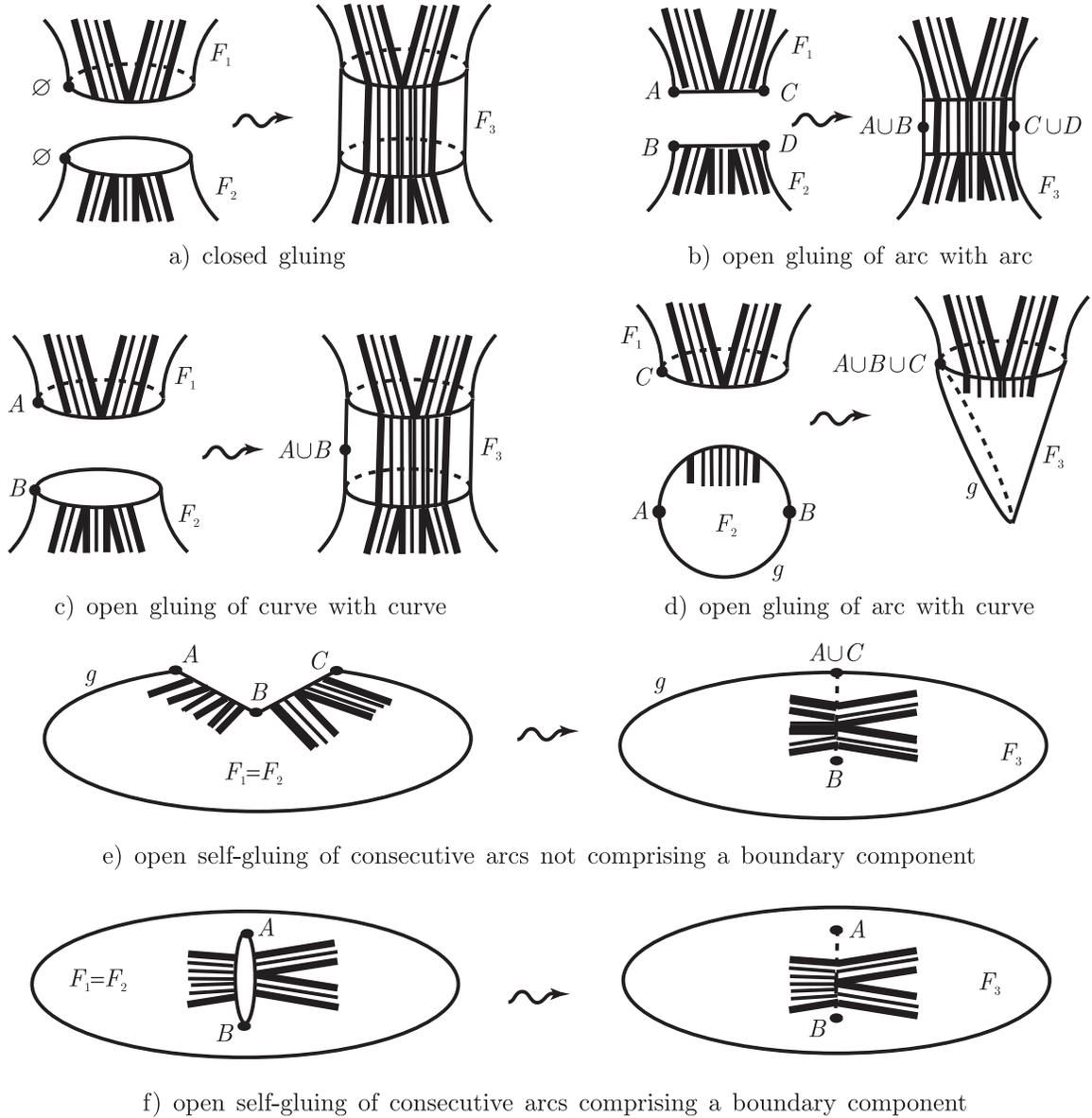}}

\caption{The dif\/ferent cases of gluing.}\label{stropfig4}
\end{figure}

\subsubsection{Closed gluing and self-gluing}
Identify the two corresponding boundary components of $F_1$ and
$F_2$, identifying also the distinguished points on them {\it and
then including this point in the resulting surface} $F_3$. $F_3$~inherits a~brane-labeling from those on $F_1$, $F_2$ in the natural
way.  We furthermore glue $\alpha _1$ and $\alpha _2$ together in
the natural way, where the two collections of foliated rectangles
in~$F_1$ and~$F_2$ which meet $w_1$ and $w_2$ have the same total
width by hypothesis and therefore glue together naturally to
provide a measured foliation ${\mathcal F}$ of a closed subsurface
of $F_3$.

\subsubsection{Open gluing}

The surfaces $F_1$ and $F_2$ are distinct, and we identify $w_1$ to $w_2$ to produce
$F_3$.
There are cases depending upon whether the closure of $w_1$ and $w_2$ is an interval or a circle.  The salient
cases are
illustrated in Fig.~\ref{stropfig4},~b--d.  In each case,  distinguished points on the boundary in~$F_1$ and~$F_2$ are
identif\/ied to produce a new distinguished boundary point in $F_3$, and the brane labels are combined, as is also
illustrated.  As before, since the
$\alpha _1$-weight on $w_1$ agrees with the $\alpha _2$-weight on $w_2$, the foliated rectangles again combine to
provide a measured foliation~${\mathcal F}$ of a closed subsurface of~$F_3$.

 {\bf Open self-gluing.}~There are again cases depending upon whether the closure of $w_1$ or $w_2$ is a~circle or an interval, but there is a further case as well when the two intervals lie in a common boundary component
and
are consecutive.  Other than this last case, the construction is identical to those illustrated
in Fig.~\ref{stropfig4},~b--d.  In case the two windows are consecutive along a common boundary component, again they are identif\/ied
so
as to produce a surface~$F_3$ with a puncture resulting from their common endpoint as in Fig.~\ref{stropfig4},~e--f, where the
puncture is
brane-labeled by the label of this point, and the foliated rectangles combine to provide a measured foliation~${\mathcal F}$ of a closed subsurface of~$F_3$.

At this stage, we have only constructed a measured foliation ${\mathcal F}$ of a closed subsurface of $F_3$, and
indeed,
${\mathcal F}$ will typically not be a weighted arc family, but the sub-foliation
${\mathcal
F}'$ comprised of leaves that meet $\partial F$ corresponds to a weighted arc family $\alpha _3$ in $F_3$.
Notice that
the $\alpha _3$-weight of any window uninvolved in the operation agrees with its $\alpha _1$- or $\alpha _2$-weight.

The assignment of $\alpha_3$ in $F_3$ to $\alpha _i$ in $F_i$, for
$i=1,2$ completes the def\/inition of the various operations.
Associativity and equivariance for bijections are immediate, and so
we have our f\/irst non-trivial example of a c/o structure; see
Appendix~\ref{appendixA} for the precise def\/inition.

\subsection{Extended gluing and the open/closed Sullivan PROP}

For the PROP case, we will need to extend the gluing to the case
where one window is active and the second window is inactive.
In this case, we glue the surfaces and distinguished points as above and
simply delete foliation in the rectangle. This may result in additional windows becoming inactive.

\begin{rmk}
We do not wish to formalize c/o PROP structures here.  A PROP gluing is given by pairing of all ``ins''
to all ``outs'' of two dif\/ferent surfaces and gluing on all of them.
 In apparent terminology, we can have the open and closed PROP substructures separate or at once. Technically there are also so-called vertical compositions which in our cases are always just disjoint unions. The brane labelling is handled in the same fashion as in the c/o case.
For more details see  Appendix~\ref{coprop}.
\end{rmk}

\begin{prop}
The open/closed Sullivan spaces are  closed under the extended
gluing, when  gluing an open $($respectively closed$)$ ``in'' to an open
$($respectively closed$)$ ``out'' window
 with the same weight or to an empty ``out'' window.
\end{prop}

\begin{cor}
 This also gives the structure of a
c/o, i/o modular operad in the terminology of the appendix.
This implies that these spaces  form c/o PROP.
These structures also exist in the brane labelled case.
\end{cor}

\begin{proof}
Both conditions for families in the Sullivan spaces are stable under
the gluing. (1)~If arcs only run from ``in'' to ``out'', they also
do so after gluing an ``in'' to an ``out'' window: Indeed a
foliation could only run from ``in'' to ``in'' on the glued surface
if there was a foliation running from ``in'' to ``in'' in the
surface to whose ``in'' window we glue. The extended gluing only
kills foliations.  (2) After gluing all ``in'' boundaries are
active: This is clear if we do not glue to an empty ``out''. But
even if we glue to an empty ``out'' this holds true, since in this
case only leaves get deleted on the surface to whose ``in'' we glue.
These foliations run to ``out'' windows of that surface and hence to
``out'' windows of the glued surface. The ``in'' windows of the
glued surface are unaf\/fected. The number of inactive ``out'' windows
may of course increase.

Since the gluing is associative, we can obtain a c/o PROP gluing by
simply gluing successively as described in the appendix. The f\/irst
gluing will be a non-self gluing, while all remaining gluings are
self-gluings. Since gluing is associative this is insensitive to
the order chosen.

The brane labelling is external to the foliation gluing, so the last statement readily follows.
\end{proof}

There is actually a open/closed colored dg-PROP structure on the
chain level, if we use cellular chains as we demonstrate in
Section~\ref{octopsection}. This is induced by a topological quasi-PROP
structure on the topological level; see the appendix for the
def\/initions of these types of PROPs.

\subsection[Discretization of the model ${\mathbb N}$-valued foliations]{Discretization of the model $\boldsymbol{{\mathbb N}}$-valued foliations}

We wish to point out that the subset of discrete valued foliations is stable under the compositions
in both the c/o structure and the c/o PROP structure.

\subsubsection[Discrete representation for discretely weighted $\beta$-arc families]{Discrete representation for {\em discretely weighted} $\boldsymbol{\beta}$-arc families}

In order to decorate, we will change the picture slightly. Previously we had arc graphs, whose
edges are not allowed to be parallel. For a {\em discretely weighted} $\beta$-arc families with weighting $wt$
we  will consider its {\em leaf representation} to be the foliation each of whose bands $e$ has
$wt(e)$ number of leaves. This means that we consider an new type of arc graph which has $wt(e)$ parallel
edges for each underlying edge of the original arc graph. We will call this the {\em discrete representative}
of the arc graph.

\section{Algebraic c/o structures}
\label{algsection}

Given a $\beta$-arc graph with weights that are natural numbers, we are going to
associate an operation on the Hochschild complexes  ${\rm CH}^*(A_{\varnothing},M_{B,B'})$
of a f\/ixed Frobenius algebra $A_{\varnothing}$ with coef\/f\/icients in a module $M_{B,B'}$ with $B,B'\in \calB$.
These modules will be given by tensor products of Frobenius algebras $A_B$ indexed by elements of $\calB$: $M_{B,B'}: A_{B'}\otimes A_B$. These operations will be def\/ined on the isomorphic
double sided bar-complexes  $B(A_B,A_{\varnothing},A_{B'})$. For most operations we will need
that $A_{\varnothing}$ is commutative, but this is not always the case.  We will indicate
when the commutativity can be dropped.

\subsection{Frobenius algebras and systems of Frobenius algebras}

\subsubsection{Notation for Frobenius algebras}
The main actors are  Frobenius algebras, so we will f\/ix some notation.

Recall that a Frobenius algebra (FA) is a triple $(A,1, \langle \; ,\;\rangle)$
where $(A,1)$ is a unital (super-) algebra and $\langle \; ,\;\rangle$ is a non-degenerate
(super-) symmetric even pairing which satisf\/ies
\begin{equation*}
\langle ab,c\rangle=\langle a,bc\rangle.
\end{equation*}

We will set
\begin{equation*}
\int a :=\langle a , 1 \rangle.
\end{equation*}
Then $\int$ is cyclically (super-)invariant, i.e.\ a trace
\begin{equation*}
\int abc = \la ab,c\ra =\la c, ab\ra=\int cab.
\end{equation*}

Since $ \langle \; ,\;\rangle$ is non-degenerate on $A$, so is $
\langle \; ,\;\rangle_{A\otimes A}:=\langle \; ,\;\rangle\otimes
\langle \; ,\;\rangle\circ \tau_{2,3}$ on $A^{\otimes 2}\otimes
A^{\otimes 2}$. Here~$\tau_{2,3}$ is the commutativity constraint
for the symmetric monoidal category applied to the second and third
factors, be it the category of vector spaces, dg-vector spaces or
$\Z/2\Z$ graded vector spaces
 In our current setup this just interchanges the second and third  factors of
$A\otimes A\otimes A\otimes A$  or in the super case changes these factors and introduces
the usual super sign.

We will omit all super signs from our discussion as they can be added in a straightforward fashion.

The multiplication $\mu:A\otimes A\to A$ has an adjoint $\Delta: A\to A\times A$
def\/ined by
\begin{equation}
\label{inveq} \langle \Delta a, b\otimes c\rangle_{A\otimes
A}=\langle a,bc\rangle.
\end{equation}

Moreover given a FA we will consider a basis $\Delta_i$, set $g_{ij}=\langle \Delta_i,\Delta_j\rangle$
and let $g^{ij}$ the coef\/f\/icients of the inverse of $(g_{ij})$, viz.\ the inverse ``metric''.

There are two special elements
\begin{equation*}
e=\mu\Delta(1)=\sum_{ij} \Delta_i g^{ij}\Delta_j,
\end{equation*}
which we call the Euler element and
\begin{equation*}
C=\Delta(1)=\sum \Delta_i g^{ij} \otimes \Delta_j,
\end{equation*}
which we call the Casimir element.

Notice that Euler element commutes with every element.
\begin{equation*}
ae=a\mu\Delta(1)=\mu\Delta(a)=\mu\Delta(1)a=ea.
\end{equation*}
This follows by direct computation and the Frobenius relations
\begin{equation*}
({\rm id}\otimes \mu)(\Delta\otimes {\rm id})=\Delta\mu=(\mu\otimes
{\rm id})({\rm id}\otimes \Delta),
\end{equation*}
which in turn follow from the def\/inition of $\Delta$ and the
invariance of the pairing~(\ref{inveq}).

\subsubsection{Adjoint maps}
Notice that for any map $r:A\to B$ between two Frobenius algebras
there is an {\em adjoint map} $r^{\dagger}:B\to A$ def\/ined by
\begin{equation*}
 \langle r^{\dagger}(b) ,a  \rangle =\langle b,r(a)\rangle.
\end{equation*}

This is equivalent to
\begin{equation}
\label{adjointeq}
 \int r^{\dagger}(b)a  =\int br(a).
\end{equation}
%In this sense $\Delta=\mu^{\dagger}$.

These maps arise
in geometric situations as follows. Let $i:N \to M $ be the inclusion map, where $M$ is a compact manifold and $N$ is a compact submanifold. Then $i$ induces a map $i^*$ from
$A:=H^*(M)$ to $B:=H^*(N)$. By Poincar\'e duality there is a push forward $i_*:B\to A$. In the previous notation if  $r=i^*$ then $r^{\dagger}=i_*$.

\begin{lem}In general, we have the Projection Formula
\begin{equation*}
r^{\dagger}(r(a)b)=ar^{\dagger}(b). % \qquad \mbox{\rm Projection Formula}
\end{equation*}
\end{lem}
\begin{proof}
\begin{gather*}
\la r^{\dagger}(r(a)b), c\ra = \la r(a)b,r(c)\ra = \la b, r(ca)\ra = \la ar^{\dagger}(b),c\ra,
\end{gather*}
where we used the cyclic symmetry of the product twice.
\end{proof}

With  the self-intersection condition Section~\ref{selfdef} in mind, we def\/ine the element
\begin{equation*}
e_r^{\perp}:=r (r^{\dagger}(1)).
\end{equation*}

\subsubsection{Systems of Frobenius algebras}
\begin{df} A $\calB$-Frobenius algebra is a set of Frobenius
algebras $A_S$ indexed by $S\in \calP(\calB)$ together with algebra maps
$r_{S,S'}:A_S\to A_S'$ whenever $S\subset S'$,
such that for $S\subset S'\subset S''$: $r_{S',S''}\circ r_{S,S'}=r_{S,S''}$.
\end{df}

Note that in particular if $A_{\varnothing}$ is commutative
every $A_S$ is an $A_{\varnothing}$ module via
the restriction map. More precisely,  every $A_S$ is a left and a right $A_{\varnothing}$
module via the maps $\lambda (a,a'):= r_{\varnothing S}(a)a'$
and $\rho(a,a') :=a' r_{\varnothing S}(a)$. If $A_{\varnothing}$ is not commutative,
we still have that $A_S$ is a left $A_{\varnothing}$ module and a right $A_{\varnothing}^{\rm op}$ module.

\subsubsection{Basic brane label systems}
Given a  brane label set $\calB$ one set choices of $\calB$-Frobenius algebras is given by
a collection~$A_B$, $B\in \calB$ and $A_{\varnothing}$ together with maps $r_B:A_{\varnothing}\to A_B$.

For any $S\in \calP\calB$ with $|S|\geq 2$ we simply set $A_S=0$ where
we allow the zero algebra to be a~Frobenius algebra.

We call these systems basic brane label systems and for simplicity deal only with these.
The data of the Frobenius algebras and morphisms will be called a basic $\calB$ Frobenius algebra.

\begin{rmk}
In the following we will mostly deal with only basic $\calB$ Frobenius algebra in order to not unduly
burden the  reader with yet more  structures that we would need in order to deal with the general case.
All the results do however generalize to the general case if we introduce propagators to f\/ix the
c/o structure on the algebraic side.
\end{rmk}

\begin{nota}
Given $A_B$ we will use the notation
$1_B$, $e_B$, $\Delta^i_B$, $\la$, $\ra_B$, $g^B_{ij}$, $g^{ij}_B$, $e^{\perp}_B$,
$\dots$ for its unit, Euler element,
basis $A_B$, metric inverse metric, $e^{\perp}_{r_B}$ etc.

We will also omit the label ${\varnothing}$ i.e.\ write  $e$ for $e_{\varnothing}$, $A$ for $A_{\varnothing}$ if no
confusion can arise.
\end{nota}

\begin{df}
\label{Eulerdef}
We say that a basic $\calB$-FA  satisf\/ies the condition of commutativity $(C)$ if $A_{\varnothing}$ is commutative.

And we say that a $\calB$-FA satisf\/ies the
the Euler compatibility condition or the {\em condition $(E)$} if for all $B\in \calB$, $a^{(1)},a^{(2)}\in A_B$
\begin{equation*}
(E) \quad \sum_{ij}\res_B^{\dagger}\big(a^{(1)}\Delta^B_i\big) g^{ij}_B
 \,  \res_B^{\dagger}\big(\Delta^B_ia^{(2)}\big)=e_{\varnothing}\res_B^{\dagger}\big(a^{(1)}a^{(2)}\big).
\end{equation*}
\end{df}

\begin{df}
\label{selfdef}
A basic $\calB$-Frobenius algebra satisf\/ies  the {\em self-intersection} condition $(I)$
if for all $B\in \calB$
\begin{gather*}
(I_1) \quad r_{B}r^{\dagger}_{B}  (a)=  a e^{\perp}_{B} \qquad \text{and}\qquad
(I_2) \quad  e_Be^{\perp}_B =  \res_B(e).
\end{gather*}
\end{df}

\begin{prop}A  basic system of  $\calB$ Frobenius algebra which satisfies the self-intersection condition $(I)$ satisfies the Euler condition $(E)$.
\end{prop}

\begin{proof} We will show that the r.h.s.\ and the l.h.s.\ of $(E)$ have
the same inner product with any element $b$ of $A$:

For all $b\in A_{\varnothing}$, $a, a'\in A_B$
\begin{gather*} \sum_{ij} \int b  \res_B^{\dagger}\big(a\Delta^{B}_i\big)
g_{B}^{ij}  \res_B^{\dagger}\big(\Delta^{B}_j a'\big)
 = \sum_{ij}g_{B}^{ij}
 \int_B \res_{B}(b)a\Delta^{B}_i
 \res_{B}\res_B^{\dagger}\big(\Delta^{B}_ja'\big) \\
\qquad{}  \stackrel{(I_1)}{=} \sum_{ij}g_{B}^{ij}
 \int_B \res_{B}(b)a\Delta^{B}_i \Delta^{B}_j
 a'e_{B}^{\perp}
 =
 \int_B \res_{B}(b)aa' e_B e_B^{\perp}\\
\qquad{} \stackrel{(I_2)}{=} \int_B \res_{B}(b)aa'\res_B(e)
 = \int b  e\res_B^{\dagger}(aa'),
\end{gather*}
where the f\/irst   equality follows from equation (\ref{adjointeq}),
the third equality from the def\/inition of $e_B$ and the fact that~$e_B$ as the Euler element commutes with all other elements and the
last equality follows from the projection formula and the fact that~$e$ commutes.
\end{proof}

\subsubsection{Geometric data}
\label{geodatasection} One example of the basic data is given by a
compact manifold $M$ together with a collection $N_B\subset M$, $B\in
\calB$ of compact submanifolds. We can then set $A_B:=H^*(N_B)$ and
use the restriction maps $r_B$     given by pullback. These satisfy
the f\/irst  condition $(I_1)$ of $(I)$ due to the self-intersection
formula where $e_B^{\perp}=e(N_{M/N_B})$ is the Euler class of the
normal bundle of $N_B$ in $M$. The second condition ($I_2$) follows
from the excess intersection formula for homology~\cite{Qu} applied
to the diagram
\begin{equation}
\label{intdiagram}
\begin{CD}
N_B @>\Delta_B>>N_B \times N_B \\
@Vi_BVV @V i_B\times i_B VV\\
M@>\Delta>> M\times M
\end{CD}
\end{equation}
keeping in mind that $\mu=\Delta^{*}$, $\res_B=i_B^*$
\begin{gather*}
\res_B(e)=\res_B\mu\mu^{\dagger}(1)=i_B^*\Delta^*\Delta_*(1)=\Delta^*_B(i_{B}^*\times i_{B}^*)\Delta_*(1) \\
\phantom{\res_B(e)}{} =e_B^{\perp}\Delta_B^*\Delta_{B*}i_B^*(1)=e_B^{\perp}\Delta_B^*\Delta_{B*}(1)=
e^{\perp}_Be_B.
\end{gather*}

Alternatively one can use decomposition $TM|_{N_B}=TN_B\oplus N_{M/NB}$ and multiplicativity of the Euler class.

\begin{rmk}
We can also formalize the geometricity by  staying in the framework of Frobenius algebras, but
postulating a new axiom which guarantees the equations one would obtain from excess intersection formulas from all embedding and intersection diagrams analogous to~(\ref{intdiagram}).
\end{rmk}

\begin{cor}
A basic geometric $\calB$ Frobenius algebra satisfies the Euler condition $(E)$.
\end{cor}

\subsubsection{Hochschild complexes}
The action on the closed sector is on the Hochschild cochain-complex of $A$.
Recall that the Hochschild chain complex of an $A$ bimodule $M$ is the complex ${\rm CH}_n(A,M)$
\[
{\rm CH}_n(A,M)=M\otimes A^{\otimes n}
\]
and whose  dif\/ferential is given by $d=\sum_i (-1)^i d_i$,
where the $d_i$ are the pre-simplicial dif\/ferentials
\begin{gather*}
d_0 (m\otimes a_1 \oto a_n ) = m a_1 \otimes\aoa{2}{n},\\
d_i (m\otimes a_1 \oto a_n) = m\otimes \aoa{1}{i-1}\otimes a_{i}a_{i+1}\otimes
\aoa{i+2}{n}, \\
d_{n} ( m\otimes a_1 \oto a_n ) = a_n m\otimes\aoa{1}{n-1}.
\end{gather*}

There are degeneracies inserting $1$ into the $i$th position
\begin{gather*}
s_i: \ {\rm CH}_n(M,A) \to {\rm CH}_{n+1}(M,A),\\
\phantom{s_i:}{} \ \ m\otimes a_1 \oto a_n \mapsto m\otimes a_1
\oto a_{i-1}\otimes  1 \otimes a_i \oto a_n.
\end{gather*}

The Hochschild chain complex ${\rm CH}_*(A,M)$ is also sometimes called
the cyclic bar complex and is denoted by $B_*(A,A)$.

The Hochschild co-chain complex is dually given by
\[
{\rm CH}^n(A,M)={\rm Hom}(A^{\otimes n},M)
\]
with the dual dif\/ferential $\d=\sum_i (-1)^i \d_i$,
$\d:{\rm CH}^n(A,M)\to {\rm CH}^{n+1}(A,M)$, where for $f\in {\rm CH}^n(A,M)$
\begin{gather*}
\d_0 f( a_1 \oto a_{n+1} ) =  a_1f(\aoa{2}{n}),\\
\d_i f(a_1 \oto a_{n+1}) = f(\aoa{1}{i-1}\otimes a_{i}a_{i+1}\otimes
\aoa{i+2}{n}), \\
\d_{n+1} ( a_1 \oto a_{n+1}) = f(\aoa{1}{n})a_{n+1}.
\end{gather*}

The degeneracies dualize to $\s_i:{\rm CH}^n(A,M)\to {\rm CH}^{n-1}(A,M)$
\begin{equation*}
\s_i f(a_1 \oto a_{n-1})=f(\aoa{1}{i-1}\otimes 1\otimes
\aoa{i}{n-1}).
\end{equation*}

In case $M=A$ multiplication of functions gives a natural product
\[
\cup: \ {\rm CH}^n(A,A)\otimes {\rm CH}^m(A,A)\to {\rm CH}^{n+m}(A,A).
\]

Notice that if $A$ is a Frobenius algebra, it is isomorphic to its
dual as a bi-module. Since $A\simeq {\check A}$, ${\rm CH}^n(A,A)\simeq
A\otimes {\check A}^{\otimes n}\simeq A^{\otimes n+1}\simeq
{\rm CH}_n(A,A)$. The cup product  and the product pairing, make
${\rm CH}^{\bullet}$ into a graded Frobenius algebra. Furthermore, the
dif\/ferentials  $d$ and $\d$ dualize to one and another.

\subsubsection{Reduced Hochschild complex}

For technical reasons discussed in \cite{hoch2} it is actually
easier to work with the reduced Hochschild co-chain complex
$\overline{\rm CH}^*(A,A)$. This complex is the subcomplex of functions
that vanish on all degeneracies that is functions $f:A^{\otimes
n}\to A$ such that $f(a_1,\dots,1,\dots, a_n)=0$ where $1$ is
plugged in into any position. The complex inherits the dif\/ferential
and computes the same cohohmology as the original complex. Dually
there is the reduced bar complex $\overline {\rm CH}_*(A,A)$ or $\bar
B_*(A,A)$ which we shall use in the closed sector.

\subsubsection{Double sided bar construction}
For the open action,
we will consider the double sided bar complexes
$B(S,T):=B_{\bullet}(A_S,A,A_T)$
whose components are def\/ined as
\[
B_n(A_S,A,A_T)=A_S\otimes A^{\otimes n}\otimes A_T
\]
and whose dif\/ferential is given by $d=\sum_i (-1)^i d_i$, where the
$d_i$ are the pre-simplicial dif\/ferentials
\begin{gather*}
d_0 (a_S\otimes a_1 \oto a_n \otimes a_T) = a_Sa_1\otimes \aoa{2}{n}\otimes a_T,\\
d_i (a_S\otimes a_1 \oto a_n \otimes a_T) = a_S\otimes
\aoa{1}{i-1}\otimes a_{i}a_{i+1}\otimes
\aoa{i+2}{n} \otimes a_T,\\
d_{n} (a_S\otimes a_1 \oto a_n \otimes
a_T) = a_S\otimes\aoa{1}{n}\otimes a_n a_T.
\end{gather*}

\begin{rmk}
Notice that since we are dealing with Frobenius algebras, this is
isomorphic to  ${\rm CH}^{\bullet}(A,A_T\otimes A_S)$ and again the
dif\/ferentials dualize to one and another.
\end{rmk}

\subsubsection{Degeneracies}
The double sided complex is actually simplicial,
which means that it also has degeneracy maps
\begin{gather*}
s_i: \ B_n(A_S,A,A_T) \to B_{n+1}(A_S,A,A_T), \\
\phantom{s_i:} \ \ a_S\otimes a_1 \oto a_n \otimes a_T \mapsto a_S\otimes a_1 \oto
a_{i-1}\otimes  1 \otimes a_i \oto a_n \otimes a_T.
\end{gather*}

\subsubsection{Brane labelled  bar complexes}

Fix a $\calB$ Frobenius algebra. For a window $w$ on a brane
labelled surface with labelling $\beta$ we set $B(w):=B(\beta (w))$
if $w$ is open and
$B(w)=B(\beta(w))=B(\varnothing,\varnothing):=\overline{\rm CH}_n(A,A)$ if $w$ is
closed.

\subsection{Gluing in brane labelled complexes}

As we have noted, the complexes $B(\beta)$ have non-degenerate
graded inner products which allows us to dualize them.
 Using this
dualization, we can compose two correlation functions for the basic
brane label case. In the general case we would need to introduce
propagators to do this. We will refrain from adding this technical
point here for clarity of the discussion.

Given  (graded) vector spaces $V_{\beta}: \beta\in I$ over a ground
f\/ield $k$ with an involution~$\bar{}$ on~$I$, isomorphisms
$V_{\beta}\simeq V_{\bar\beta}$ s.t.\ $\bar{\bar{a}}=a$ and (graded)
non-degenerate even symmetric pairings $\la\;,\;\ra_{\beta}$ for
each $V_{\beta}$
 let $C_{\beta}=\sum
\Delta^{\beta}_i g^{ij}_{\beta} \otimes \bar\Delta^{\beta}_j$ be the
Casimir element for the induced pairing between $V_{\beta}$ and
$V_{\bar\beta}$ expressed in a basis $(\Delta^{\beta}_i)$.

We can compose two correlators two correlators $Y:\bigotimes_{l \in
L} V_{\beta_l} \to k$ and $Y':\bigotimes_{l'\in L'}
V_{\beta_{l'}}\to k $ where $L$ and $L'$ are labelling sets by
inserting a Casimir.

More precisely given $l_0\in L$ and $l'_0\in L'$ such that
$\beta_{l_0}=\bar\beta_{l_0'}=\beta$ we def\/ine their composition
$Y\circ_{l_0l_0'} Y': \bigotimes_{l''\in (L \setminus \{l_0\} \sqcup
L' \setminus \{l_0'\})} \to k$ by
\begin{gather*}
Y\circ_{l_0l_0'} Y' \Big(\bigotimes_{l\in L\setminus \{l_0\}} v_l\otimes
\bigotimes_{l'\in L'\setminus \{l_0'\}} v_{l'}\Big) \\
\qquad{}
=\sum_{ij}
Y\bigg(\Big(\bigotimes_{l\in L\setminus \{l_0\}} v_l\Big)\otimes
\Delta^{\beta}_i\bigg)g^{ij}Y'\bigg(\Big(\bigotimes_{l'\in L'\setminus
\{l_0'\}} v_{l'}\Big)\otimes \bar \Delta^{\beta}_j\bigg),
\end{gather*}
where $\Delta^{\beta}_i$ is in the position $l_0$ and $\bar
\Delta^{\beta}_j$ is in position $l_0'$. If $\beta_{l_0}\neq
\bar\beta_{l'_0}$ we set the composition to $0$. In the more general
case these would be non-zero and the composition would use
propagators.

Alternatively, we could also dualize the maps $Y$ using
$C_{\beta_l}$ in any position $l$ to obtain maps to~$V_{\beta_l}$
instead of $k$ and then compose these maps. There is an analogous
procedure for self-gluing.

\begin{rmk}
In our case $I=\calB\times \calB\cup \{\varnothing,\varnothing\}$,
$\overline{(S,T)}=(T,S)$, $V_{\beta}=B(\beta)$ the corresponding bar
complex and the isomorphisms $\bar{}:B(S,T)\to B(T,S)$  are given by
\begin{equation*}
a_S\otimes \aoa{1}{n}\otimes a_T\mapsto a_T \otimes
\aoa{n}{1}\otimes a_S
\end{equation*}
and $B(\varnothing,\varnothing)\to B(\varnothing,\varnothing)$
\begin{equation*}
 \aoa{0}{n} \mapsto a_0 \otimes
\aoa{n}{1}.
\end{equation*}
\end{rmk}

\begin{rmk}
In order to give a brane labelled c/o structure, we should consider
slightly enriched more complicated data. For this we would have to
look at tensors products of cyclic tensor products of bar complexes
$B(S_1,S_2)\otimes B(S_2,S_3)\oto B(S_n,S_1)$. Again in the interest
of brevity, we will not introduce this kind of complexity in a
formal fashion here.
\end{rmk}

\begin{rmk}
There are actually two c/o structures, one can compose the bar
complexes or their duals, viz.\ the correlators. Of course these
operations are dual to each other.
\end{rmk}

\section{Correlators}
\label{corsection}

\subsection{A universal formula for  correlators}
 There is a universal formula for the correlators. It is given by partitioning, decorating
and  de\-com\-posing the surface of a discretely weighted arc family
  along the arcs into little pieces of surface $S_i$ and integrating around
these pieces. We will now give the details.

\subsubsection{Decorating the boundary}
Fix a basic $\calB$ collection of Frobenius algebras.  For each {\em
discretely weighted} $\beta$-arc family $\alpha$,  we will def\/ine a
map
\[
Y(\alpha):= \bigotimes_{w \in \mathrm {Windows~of~}\alpha}B(\beta(w))\to k.
\]
These functions are homogeneous and their homogeneous components  are  zero by def\/inition unless
 $a_w\in B_{\alpha(w)-1}(\beta(w)) $.

Given a collection of homogeneous elements
we will decorate the pieces belonging to the boundary of the discrete representation of $\alpha$ by
the elements of the bar complexes
\begin{gather*}
a_w = a^w_S\otimes a_1^w \oto a^w_{\alpha(w)-1}\otimes a_T^{\prime w} \qquad \text{if} \quad \beta(w)=(S,T),\\
a_w = a^w_0\otimes a_1^w \oto a^w_{\alpha(w)-1}  \qquad \text{if} \quad
\beta(w)=(\varnothing,\varnothing).
\end{gather*}

Notice that  the boundary of the underlying surface minus the
discrete representative of the graph is a disjoint union of
intervals, which may or may not contain marked points. We call these
the {\em boundary pieces}. There are three types of boundary pieces
\begin{enumerate}\itemsep=0pt
\item[(1)] those not containing a marked point,
\item[(2)]  those
containing a marked point $\beta$-labelled by $\varnothing$,
\item[(3)] those containing a marked point with
$\beta$-label not $\varnothing$.
\end{enumerate}
In case~(3), if we remove the marked point we will have two components
which we will call half sides of the boundary piece. Now
each piece of type~(1) and the  half sides of the pieces of type~(3)
 belong to a unique window.
A piece of type~(2) comes from a unique closed window/boundary component. Moreover these
pieces all come in a natural linear order in each window as do the half sides of a piece of type (3) if
we consider the marked point to lie in between the half sides.

 We decorate
the boundary pieces as follows:

\vspace{1mm}

{\centering \begin{tabular}{@{\extracolsep{0pt}}c|l|l}
Type&description of $s$ &Decoration\bsep{1pt}\\
\hline
 (1) &$s$ is the $i$th piece of type (1) of the window $w$&$a_i^w\in A$\tsep{1pt}\bsep{1pt} \\
(2)& $s$ is the unique piece of type (2) of the closed window $w$&$a_0^w$ \bsep{1pt} \\
 (3)&the marked point
of  $s$ is labelled by $S$ and the half &
$(a_S^{\prime w_1},a_S^{w_2})$\\
&sides in their order belong to the &\\
&not necessarily distinct
 windows $w_1$, $w_2$&
\end{tabular}

}

%maybe A_empty=A or A_emtpy-C

\subsubsection{Weights}

 Let $\{S_i:i\in I\}$ be the components of
complement of the discrete representative of a given discretely weighted
$\beta$ arc family $\alpha$.
Each of these pieces has a polygonal boundary, where the sides of the polygons alternate between pieces
of the boundary and arcs running
between them. If we decorate the surface as described above every second side is a decorated
piece of boundary. We will call these the decorated sides.

To each decorated side $s$ of $S_i$, we associate a weight depending $\wt$
on its type and decoration.

\vspace{-2mm}

\begin{table}[h]\centering
\caption{General weights.}
\vspace{1mm}
\begin{tabular}{lll}
Type&Decoration& Weight $\wt(s)$\bsep{1pt}\\
\hline
(1) $s$ without marked point&$a\in A_{\varnothing}$&$a$\tsep{1pt}\bsep{1pt}\\
(2) $s$ with marked point marked by $\varnothing$&$a\in A_{\varnothing}$&$a$\bsep{1pt}\\
(3) $s$ marked point marked by $S$&$(a_S^{(1)},a_S^{(2)})$, $a_S^i\in A_{S}$&
$r^{\dagger}_{\varnothing S}(a_S^{(1)} a_S^{(2)})$
\end{tabular}
\end{table}

\subsubsection{The formula}
Moreover, the $S_i$ are oriented and so hence are their boundaries.
This means that the sides comprising each boundary component come with a cyclic order.
If there is only one boundary component for a given $S_i$, this gives a cyclic order over which we will integrate the given weights.
If there are more components, in which case the underlying arc family is not quasi-f\/illing,
then we need to assume (C) in order to make the following expression independent of choices.
For a homogeneous ${\mathbf  a}= \bigotimes_{w \in \mathrm {Windows~of~}\alpha} a_w \in \bigotimes_{w \in \mathrm {Windows~of~}\alpha}B(\beta(w))$ such that  $a_w\in B_{\alpha(w)-1}(\beta(w))$,
we decorate as above and def\/ine
\begin{equation}
\label{cordef}
Y_{S_i}({\mathbf a})=\int e^{-\chi(S)+1}\prod_{\substack{\text{Decorated sides}\\
\text{$s$ of $S_i$}}} \wt(s)  \prod_{\substack{\text{Punctures $p$}\\ \text{inside $S_i$}}}
r^{\dagger}_{\varnothing \beta(p)}(e_{\beta(p)}).
\end{equation}
If $\mathbf a$ is as above but there is  some $a_w\notin B_{\alpha(w)-1}(\beta(w))$ we set $Y_{S_i}(\mathbf a)=0$.

We then def\/ine{\samepage
\begin{equation}
\label{cordefgraph}
Y_{(\Gamma,w))}({\mathbf a}):=\prod_i Y_{S_i}({\mathbf a})
\end{equation}
and extend by linearity.}

\subsubsection{Signs} The correlators above actually have hidden signs which come
from the permutation of the input variables to their respective position. In the bar complexes these signs can be read of\/f by imposing that the tensor symbols have degree~1. On the geometric side
there are signs as well which are f\/ixed by f\/ixing an enumeration of the f\/lags, angles or edges.
In general we adhere to the sign conventions spelled out in \cite[Section 1.3.4]{hoch2}.

\subsection{Action of the c/o structure of discretely weighted arc-graphs}

We say a c/o structure acts via correlation functions if the
composition of the elements of the c/o is compatible with the
composition of the correlation functions. In short $Y(\alpha
\circ_{w,w'} \alpha')=Y(\alpha)  \circ_{w,w'} Y(\alpha')$ where
$\circ_{w,w'} $ denotes the gluing of the window $w$ and $w'$ holds
as well as the corresponding equation for the self-gluing.

\begin{thm}
\label{discretethm} For a basic  $\B$-Frobenius algebra the c/o
structure of discretely weighted arc-graphs acts  on  the
collection of complexes $B(\beta)$ and the isomorphic Hochschild
complexes via the correlation functions $Y$.
\end{thm}

Just like there are algebras over operads, we can def\/ine algebras
over c/o structures.  The theorem above reads: The collection of bar
complexes $B(S,T)$ form an algebra over the c/o structure of
discretely weighted arc-graphs.

\begin{proof}
The proof is a case by case study which occupies Section~\ref{casessection}.
\end{proof}

\subsection{Case by case analysis  of the discrete arc-graph action}
\label{casessection}

\begin{figure}[th]
\centerline{\includegraphics[scale=0.96]{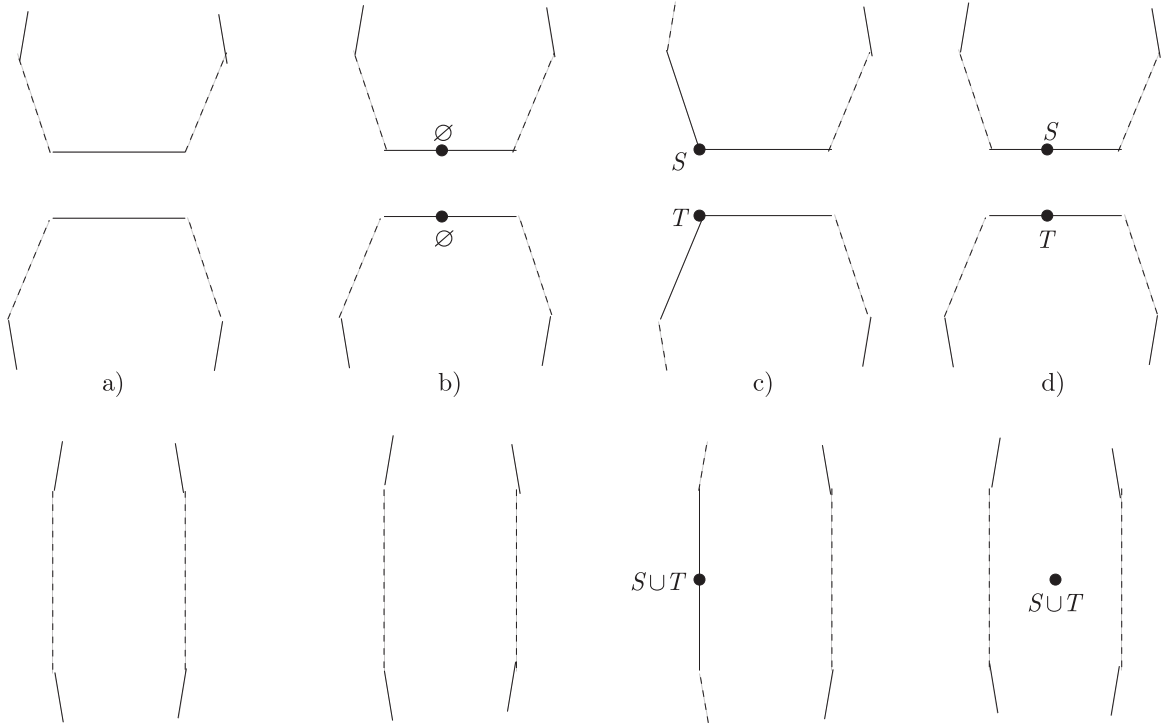}}
\caption{Types of local gluings: a)~two sides
without marked points; b)~side with marked point labelled by $\varnothing$
to side with marked point labelled by $\varnothing$; c)~half a labelled
side with marked point labelled by $S$ to half side with marked
point labelled by $T$; d)~full side with labelled point marked by $S$
to full side with labelled point marked by $T$.}\label{gluings}
\end{figure}

\begin{figure}[th]
\centerline{\includegraphics[scale=0.96]{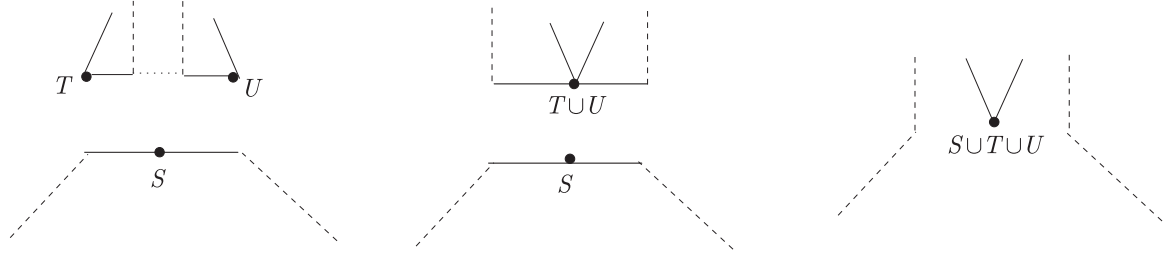}}
\caption{e)~Gluing of a lone window with an open
marked point to window with two marked points.}\label{gluingstwo}
\end{figure}

The gluing of the surfaces with discrete arcs breaks down into
individual local gluings of pairs of surfaces $S_i$ and $S'_j$. In
case the surfaces are distinct, there are f\/ive cases of this gluing
depicted in Figs.~\ref{gluings} and~\ref{gluingstwo}. The
cases~a) and~b) are the ones familiar from the closed gluing~\cite{hoch2}, the cases~c) and~d) are new in the open/closed case.
The gluing~c) appears when we are gluing two open windows where none
of them is the only window in its boundary component. The gluing~d)
appears when gluing two open windows each of which is the only
window in its boundary component. The most complicated case is when
one of the windows is the only window, while the other is not. This
case~e) is given in Fig.~\ref{gluingstwo}.

\subsubsection{Non-self gluing, the closed case}
In the case that there is no self-gluing: a)--b) correspond to:
\begin{equation}
\bigg[\int \dots w  \Delta^{\varnothing}_i  w' \dots \bigg] g^{ij}_{\varnothing}
\bigg[\int  \dots w''  \Delta_j^{\varnothing} w''' \dots\bigg]
=\int w' \dots w'' \dots w''',
\end{equation}
where the $w$, $w'$, $w''$ are the weights of the adjacent sides and we use Einstein summation conventions. We have also included
any factors of~$e$ into the ellipses.

\subsubsection{Non-self gluing case in the simple brane case}

Assuming the simple brane case, all labels have to be the same, say
$S$ then c) corresponds to
\begin{gather*}
\bigg[\int \dots w  \idagger_{S}(a_S \Delta^{S}_i  w'\dots\bigg]  g^{ij}_{S}
\bigg[\int  \dots w'' \idagger_S (\Delta_j^{S}a'_S) w''' \dots\bigg]  \\
\qquad{} = \bigg[\int_S \res _S (w' \dots w) a_S \Delta^{S}_i \bigg] g^{ij}_{S}
\bigg[\int_S  \Delta_j^{S}a'_S \res_S( w''' \dots w'')\bigg]\\
\qquad{} = \int w' \dots w \idagger_S (a_S a'_S)  w''' \dots w''
 =  \int \dots w \idagger_S(a_S a'_S)  w'''  \dots w'' w',
\end{gather*}
while d) corresponds to
\begin{gather*}
\bigg[\int \dots w  \idagger_{S}(\Delta^{S}_i  \Delta^{S}_{k}) w'\dots \bigg] g^{ij}_{S}g^{kl}_{S}
 \bigg[\int  \dots w'' \idagger_{S}( \Delta_l^{S} \Delta_j^{S}) w''' \dots \bigg]  \\
\qquad{} = \int w' \dots w \idagger_S (\Delta_i^S g_S^{ij}\Delta_j^S)  w''' \dots w''
 = \int \dots w \idagger_S(e_S)  w'''  \dots w'' w'.
\end{gather*}

The case e) has two subcases. 1) there are three surfaces which are
glued:
\begin{gather*}
\bigg[\!\int \dots w  \idagger_{S}(a_S \Delta^{S}_{i}) w'\dots \bigg]
g^{ij}_{S}
\bigg[\!\int \dots w''  \idagger_{S}( \Delta^{S}_{k} a'_S) w'''\dots \bigg] g^{kl}_{S}
\bigg[\!\int  \dots w^{\rm (vi)} \idagger_{S}( \Delta_j^{S} \Delta_l^{S}) w^{\rm (v)} \dots \bigg]   \\
\qquad{}\times\bigg[ \int_S a'_S\res_S( w'''\dots w'' )\Delta^{S}_{k} \bigg ] g^{kl}_{S}
\int _S  \Delta_l^S \res_S\big(w^{\rm (v)} \dots w^{\rm (vi)}\big)\Delta_j^{S}
  g^{ij}_{S}\bigg[\int_S  \Delta^{S}_{i}\res_S( w' \dots w)a_S\bigg]   \\
\qquad{} = \int_S \dots w\idagger_S(a_Sa'_S) w'''\dots w'' w^{\rm (v)} \dots w^{\rm (vi)} w',
\end{gather*}
and the case 2) where only two surfaces are glued %better words
\begin{gather*}
\bigg[\int \dots w  \idagger_{S}(a_S \Delta^{S}_{i}) w' \dots w''
\idagger_{S}\big( \Delta^{S}_{k} a'_S\big) w'''\dots \bigg] g^{ij}_{S}
g^{kl}_{S} \bigg[\int  \dots w^{\rm (vi)} \idagger_{S}\big( \Delta_j^{S} \Delta_l^{S}\big) w^{\rm (v)} \dots \bigg]   \\
\qquad{} =\bigg[
 g^{ij}_{S} \int_S a'_S\res_S\big( w'''\dots w  \idagger_{S}\big(a_S \Delta^{S}_{i}\big) w' \dots w'' \big)\Delta^{S}_{k}  \bigg] g^{kl}_{S}
 \int _S \Delta_l^S \res_S\big(w^{\rm (v)} \dots w^{\rm (vi)} \big) \Delta_j^{S}
 \\
\qquad{} = g^{ij}_{S} \int_S a'_S\res_S\big( w'''\dots w  \idagger_{S}\big(a_S \Delta^{S}_{i}\big) w' \dots w'' w^{\rm (v)} \dots w^{\rm (vi)} \big) \Delta_j^{S} \\
\qquad{}=  g^{ij}_{S}\int  \dots w  \idagger_{S}\big(a_S \Delta^{S}_{i}\big) w' \dots w'' w^{\rm (v)} \dots w^{\rm (vi)} \idagger\big(\Delta_j^S a_S'\big)w''' \dots.
\end{gather*}

In this case we have to use commutativity $(C)$ and
the Euler compatibility $(E)$:
\begin{gather*}
  g^{ij}_{S}\int  \dots w  \idagger_{S}\big(a_S \Delta^{S}_{i}\big) w' \dots w'' w^{\rm (v)} \dots w^{\rm (vi)} \idagger_S\big(\Delta_j^S a_S'\big)w''' \dots \\
\qquad {}\times g^{ij}_{S}\int  \dots w   w' \dots w'' w^{\rm (v)} \dots w^{\rm (vi)}\idagger_{S}\big(a_S \Delta^{S}_{i}\big) \idagger_S\big(\Delta_j^S a_S'\big)w''' \dots \\
\qquad{} =\int \dots w  w' \dots w'' w^{\rm (v)} \dots w^{\rm (vi)} \idagger_S(a_Sa'_S) e,
\end{gather*}
which is the contribution we get from the glued surfaces, since there is one self-gluing involved  and
this makes the Euler characteristic go up by one.

\subsubsection{The self gluing cases}
\label{twocasespar} So far we have assumed that the two (half) sides
that are glued are on dif\/ferent $S_i$ and $S'_j$. It can happen that
they belong to the same surface. In these cases, much like in the
case e) 2), there are fewer integrals and instead an Euler class
factor $e$ appears. In the case that there is self-gluing: a)--b)
correspond to:
\begin{equation*}
\bigg[g^{ij}_{\varnothing}  \int \dots w  \Delta^{\varnothing}_i  w'    \dots w''  \Delta_j^{\varnothing} w''' \dots\bigg]
=\int w' \dots w'' \dots w''' e.
\end{equation*}
This r.h.s.\ is the contribution to  correlator for the glued
surface since the self-gluing
 changes the Euler characteristic by $-1$. Again we need to use $(C)$.

 The cases c), d) are analogous to the case e) 2). The self gluing decreases the Euler characteristic by one,
 while the summation gives the factor $e$ by condition~$(E)$. The two cases for e) then either involve only one or two integrals, respectively; correspondingly the gluing then gives rise to a factor of $e$ or $e^2$,
respectively.

There is one more local gluing which comes from gluing consecutive windows corresponding
to Fig.~\ref{stropfig4} cases e) and f). In these cases two half sides of a {\em single}
side marked by a point labelled by some $S\neq \varnothing$ are glued together.

\begin{figure}[th]
\centerline{\includegraphics[scale=0.97]{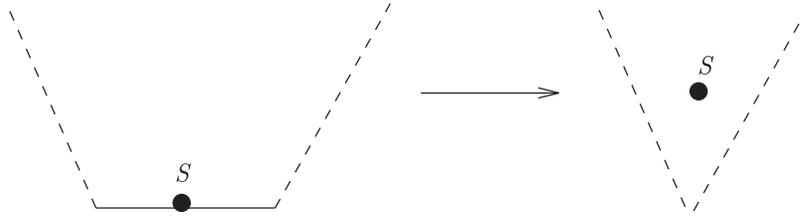}}
\caption{(g) Gluing consecutive half sides.}\label{gluingsthree}
\end{figure}

This corresponds to the following equation in the contribution to the correlator of $S_i$
\begin{equation*}
\int \dots w \idagger\big(\Delta_i^Sg^{ij}\Delta_j^S\big)w' \dots=
\int \dots w \idagger(e_S) w'' \dots.
\end{equation*}

When performing the gluing on the whole window, it might happen,
that there are also self-gluings for the surfaces $S_i$ at some
other sides, in which case we need $(C)$ and $(E)$ and proceed as above.
If there are more consecutive open gluings on one $S_i$, we produce
two punctures and correspondingly two factors of~$\idagger_S(e_S)$.

\subsection{Partitioning arc graphs and the action of arc graphs}
\label{arcgraphactionpar}
Given a discrete  weighting $w$ for an arc graph $(F,\beta,\Gamma)$ or $\Gamma$ for short
we def\/ine  $\Gamma(w)$ to be the graph in which the edge $e$ has been duplicated $w(e)-1$ times.
This is we replace $e$ with $w(e)$ parallel copies of $e$.

We def\/ine  its discretized version $P\Gamma$ of $\Gamma$
to be given by the  formal sum of the $\Gamma(w)$
\begin{equation*}
P\Gamma=\sum_{w: \ {\rm discrete~weighting~of~}\Gamma} \Gamma(w).
\end{equation*}

\subsubsection{Action of arc graphs}
\label{graphactionpar}
Given correlation functions $Y(\Gamma,w)$ for discretely weighted arc graphs $(\Gamma,w)$
we def\/ine
\begin{equation*}
Y(\Gamma):=Y_{P\Gamma}
\end{equation*}
by extending $Y$ as a function to formal sums.

\subsubsection{Examples: multiplication and comultiplication in the  open sector}
\begin{figure}[t]
\centerline{\includegraphics[scale=0.96]{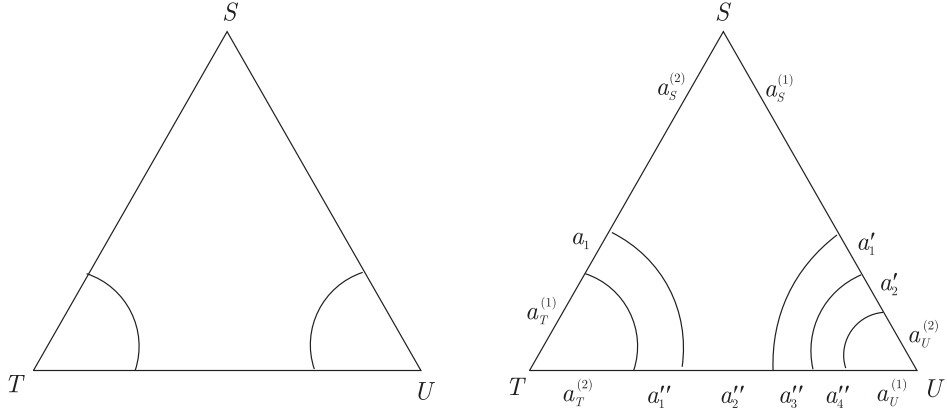}}
\caption{To the left: an arc graph $\Gamma$ in a triangle yielding the open sector multiplication or  the co-multiplication. To the right: a  discrete summand of $P\Gamma$ with weights $2$ and $3$ together with a decoration of it.}\label{openmult}\vspace{-2mm}
\end{figure}

As an example we will consider the arc graphs given in Fig.~\ref{openmult}. This triangle gives a correlation function, which
when dualized on the bottom edge yields a multiplication and when
dualized on the top two edges yields a comultiplication. A more
familiar form is given in the Sullivan case; see Section~\ref{sullexpar} below.

Given
\[
a=a^{(2)}_T\otimes a_1\oto a_n\otimes a^{(1)}_S\in B_n(A_T,A,A_S)
\]
and
\[
a'=a_S^{(2)}\otimes a'_1\oto a'_m\otimes a^{(1)}_U\in
B_m(A_S,A,A_U)
\] their product which lies in $B_{n+m+1}(A_T,A,A_U)$
is given by
\begin{equation*}
m_{TSU}(aa') = a_T^{(2)}\otimes a_1\oto a_n\otimes   r^{\dagger}_{S}(a^{(1)}_Sa^{(2)}_S) \otimes a'_1\oto a'_m\otimes a^{(1)}_U.
\end{equation*}
The calculation goes as follows. The integrals are $\int
r^{\dagger}_T (a^{(1)}_T a^{(2)}_T)=\int_T a^{(1)}_T a^{(2)}_T$
which dualizes to ${\rm id}_T(a_T^{(1)})=a_T^{(1)}$ and likewise for $U$.
For the rectangles we get either $\int a_i a''_j$ or $\int a_i'
a''_k$ which dualize to ${\rm id}(a_i)$ and ${\rm id}(a_j)$. Finally we get
$\int r_S^{\dagger}(a_S^{(1)}a_S^{(2)})a''_{n+1}$ which dualizes to
$r_S^{\dagger}\circ \mu$.

The corresponding  co-products  by dualizing are $\Delta_{TSU}: B(A_T,A,A_U)\to B(A_U,A,A_S) \otimes B(A_S,A,A_T)$ let \begin{gather*}
a''=a_T^{(2)}\!\otimes a''_1\oto a''_{n+m+1} \otimes a_U^{(1)},\\
\Delta_{TSU}(a'')=
\sum_i\! \big[a_T^{(2)}\!\otimes a''_1\oto a''_{i-1}\!\otimes (r_S(a_i))^{(1)}\big]
%\\\phantom{\Delta_{TSU}(a'')=}{}
\!\otimes\! \big[(r_S(a_j))^{(2)}\! \otimes a_{i+1}
\oto a_n\!\otimes a^{(1)}_U \big],
\end{gather*}
where $\Delta_S(r_S(a_i)=(r_S(a_i))^{(1)})\otimes (r_C(a_j))^{(2)})$ using Sweedler's notation.
For this we dualize $\int r_S^{\dagger}(a_S^{(1)}a_S^{(2)})a''_{n+1}=\int_{S\otimes S}(a_S^{(1)}\otimes a_S^{(2)} \Delta_S(r_S(a''_i))$.

\section{Open/closed string topology}
\label{octopsection}
\subsection{Cell actions}
We will now consider chain complexes whose chain groups are generated by generators indexed
by arc graphs. The basic way to obtain actions of these chain complexes is to
let each generator act via the graph that indexes it  as we explain below.

Analogously to the purely closed situation \cite{hoch1,hoch2} there are two cases which we can study.
The open/closed moduli space case and the open/closed Sullivan PROP case. Although their basic underpinnings
are the same the details are slightly dif\/ferent, again as in~\cite{hoch1,hoch2}. For the open/closed PROP case,
we will have to change the actions of the weighted arc graphs slightly by using degeneracy maps.
This is what corresponds to the breaking of the symmetry between ``ins''
and ``outs'' in the PROP itself. After putting in these degeneracies, we
obtain a dg-PROP action on the cell level. We will give the full details below. The moduli space case is discussed further in Section~\ref{modulisection}.

\subsection{The open/closed Sullivan  c/o-colored topological quasi-PROP}
So far we have only a partial gluing structure. At the expense of
having associativity only up to homotopy, we can rectify this
partial structure to a full open/closed colored structure (see the
Appendix for a def\/inition). Again without being too technical this
means that we will have PROP gluings for two given elements  $\alpha
\in \wSull(n_1,n_2,m_1,m_2)$ and $\alpha'
\in\wSull(n_2,n_3,m_3,m_4)$ and a~paring of the closed ``out''
windows of $\alpha$ and the ``in'' windows of $\alpha$. Likewise
given elements and appropriate pairings there are gluings $\alpha
\in \wSull\!(n_1,n_2,m_1,m_2)$ and $\alpha'
\in\wSull\!(n_3,n_4,m_2,m_3)$ on all the open windows. Or even gluing
all ``in'' to all ``out'' windows $\alpha
\in\wSull(n_1,n_2,m_1,m_2)$ and $\alpha'
\in\wSull(n_2,n_3,m_2,m_3)$.

This will be done by means of a f\/low. This f\/low analogous to the f\/low in \cite{hoch1}, but
dif\/fers from the f\/low in~\cite{KP}. The f\/low in~\cite{KP} might take us outside the Sullivan spaces.

Given two elements $\alpha \in \wSull(n_1,n_2,m_1,m_2)$ and $\alpha'
\in \wSull(n_3,n_4,  m_2, m_3)$, the f\/low depends on the
choice of a pairing of open ``in'' windows of $\a'$ and ``out''
windows of $\a$ and  scales the weights of all the arcs incident to
the $m_2$ open ``in'' windows of $\alpha'$ simultaneously. Given
such a pairing the f\/low at time $t$ scales each weight of an arc to
an ``in''  window $w'$ of $\a'$ by the factor of
$1-t(\alpha(w)/\alpha(w')-1)$  where $w$ is the window of $\a$
paired with $w$. At time~1, each window has the weight of its
partner under the paring.  Now glue using the previously established
gluings on all windows. Since the partial structure was bi-modular,
it does not matter in which order the gluings are performed. We can
repeat the analogous procedure for open windows or all windows at
once. The proof that these gluings are associative up to homotopy,
which is the def\/inition of a topological quasi-PROP goes along the
same line of arguments as in~\cite[Section~5.6]{hoch1}:

\begin{prop}
The open/closed Sullivan spaces $\wSull(n_1,n_2,m_1,m_2)$
form a $($two-colored$)$ brane labelled topological quasi-PROP.
\end{prop}

\begin{proof}
The two colors are open and closed. We can either choose to glue only these, or glue both open
and closed windows at one. The brane labelling is just given by the left and right brane labels of
each window. The associativity up to homotopy comes from the fact that we used a f\/low. Flowing
backwards interpolates between the dif\/ferent bracketings. These f\/low of course is only
on the non-deleted arcs, which are the same set in both bracketings. An arc is deleted if it passes
through the preimage of the glued windows. This condition is the same for both iterations
in the associativity check.
\end{proof}

Notice that we have rectif\/ied the partial structure on the topological level, but we had to pay the price
of relaxing associativity.  This weaker structure
of course induces a strict structure on the homology.

\begin{cor}
The homology open/closed Sullivan spaces $\wSull(n_1,n_2,m_1,m_2)$
form a $($two-colored$)$ brane labelled PROP.
\end{cor}

The surprising fact is that although there is only the weaker structure on the topological level,
there is already a strict structure on the chain level, when using the correct chains. This is the underlying principle of our constructions.

\subsection[A CW model for $\Sull$]{A CW model for $\boldsymbol{\Sull}$}

This paragraph is an application of the methods set forth in \cite{hoch1}.
We def\/ine the following subspaces of $\wSull(n_1,n_2,m_1,m_2)$  we let
${\Sull}_1(n_1,n_2,m_1,m_2)$ be the subspace of all $\alpha\in\wSull(n_1,n_2,m_1,m_2)$
such that the {\em $\alpha$ weight of each ``in'' window is $1$}.

\begin{prop}
${\Sull}_1(n_1,n_2,m_1,m_2)$ is a CW complex, is a sub-topological quasi-PROP and is a deformation retract of $\wSull(n_1,n_2,m_1,m_2)$.
The dimension $k$ cells of this complex are indexed by arc graphs of Sullivan type with $k+1$ arcs
and their attaching maps are given by deleting arcs and identifying this boundary with the cell of lower
dimension.
\end{prop}

\begin{proof}
The fact that this is a deformation retract is shown again by using a f\/low. This time the f\/low
scales the weights of the arcs incident to an ``in'' window  $w$ by a factor of $1-t(1/\alpha(w)-1)$.

Since for a graph of Sullivan type, all the arcs can be enumerated
by going along the ``in'' windows, we see that before taking PMC
orbits, the graphs of a given Sullivan type with $\alpha$ weight on
the ``in'' windows is simply a product of simplices. Now there is no
PMC isotropy inside this product of simplices, but under the action
sides may become identif\/ied. So for each arc graph $\Gamma$  of
Sullivan type on a brane labelled windowed surface $F$, we obtain a
cell $C(\Gamma)$ whose interior $\dot C(\Gamma)$ is given by a
product of open simplices $\dot\Delta^k$:
\begin{equation*}
\dot  C(\Gamma) =\prod_{w \in \{\text{ ``in'' windows of $F$\}}} \dot\Delta^{|\{\text{arcs incident to $w$}\}|}.
\end{equation*}
The statement about the attaching maps follows directly from the topology in $\wSull$,
where we simply delete
an arc in the limit where its weight goes to zero.
The fact that this is a sub-topological quasi PROP is immediate upon
noticing that the property that the $\alpha$ weight of each ``in'' window is  one is stable under the operation of gluing.

Since the gluing a window with $n$ arcs to a window with $m$ arcs produced at most $n+m-1$ arcs
we see that the gluing maps are indeed cellular and there are induced maps on the cellular level.

It remains to prove that these maps are associative on the nose. This follows in the same way as in \cite[Theorem~5.33]{hoch1}.
The proof there essentially goes over to the current situation for the closed part. The open part is actually simpler,
since we do not have to worry about the condition of ``twisted at the boundary'', since we keep the punctures
upon gluing.
\end{proof}

\begin{thm}
The cellular chains of $\Sull_1(n_1,n_2,m_1,m_2)$ are an open/closed
colored brane labelled PROP cell model for the spaces
$\wSull(n_1,n_2,m_1,m_2)$.
\end{thm}

\begin{proof}
This follows directly from the proposition.
\end{proof}

\begin{cor}
There is a open/closed brane labelled PROP structure on the free Abelian group generated by arc graphs of Sullivan type induced by the corresponding structure on the cellular chains of $\Sull_1$.
\end{cor}

\begin{proof}
This is def\/ined simply by the identif\/ication of the free Abelian groups of cellular chains which are generated
by the $C(\Gamma)$ and the free Abelian groups generated by the respective graphs.
\end{proof}

Notice that in this structure when gluing $\Gamma$ and $\Gamma'$ we obtain all the graphs
that can appear combinatorially by giving arbitrary weights in $\Gamma$ and $\Gamma'$ matching
on the windows that are glued, with the extra condition that the arcs of the glued graph have
the maximal number, i.e.\ the corresponding cell has the maximum possible dimension.

\subsection{The open/closed string topology action\\ on brane labelled Hochschild complexes}

\subsubsection{The correlators in the Sullivan case}
In the Sullivan graph case, when decorating and calculating the weights, we also
distinguish between ``in'' and ``out'' boundaries.
We will use the  correlators $Y^{\rm i/o}(\Gamma)$
which are obtained from the $Y((\Gamma,w))$ by using the degeneracies.

Given an discretely weighted arc family of Sullivan type
 $(\Gamma, w)$ on a surface $F$ with ``in'' and ``out'' boundary  markings.
Let
\[
\ba=\ba^{\rm (in)}\otimes \ba^{\rm (out)}\in \bigotimes_{w \in \text {``in'' windows~of~}\alpha}B(\beta(w))\otimes \bigotimes_{w \in  \text {``out'' windows~of~}\alpha}B(\beta(w)),
  \]
  we def\/ine $s\ba$
as $\ba^{\rm (in)}\otimes s\ba^{\rm (out)}$ where $s\ba^{\rm (out)}$ is def\/ined as
 follows:

Let $\ba^{\rm (out)}=\bigotimes a_w$. On the window $w$ of the discrete
family enumerate all  components
of $w\setminus\{\text{endpoints of arcs incident to  $w$}\}$ in the order
induced by the orientation of the surface starting at~$0$. Let $n_1<n_2 \dots n_k$ be the
 boundary pieces between non-parallel arcs, which are not f\/lags. Then
 we set
\begin{equation*}
sa_w :=s_{n_1}s_{n_2}\dots s_{n_k}a_w.
\end{equation*}

Given $\ba\in \ba^{\rm (in)}\otimes \ba^{\rm (out)}$ we def\/ine
\begin{equation*}
Y^{\rm i/o}(\ba):=Y(s\ba).
\end{equation*}

 According to this the correlators
will again be multilinear maps
\begin{equation*}
Y^{\rm i/o}(\alpha):= \bigotimes_{w \in \mathrm {Windows~of~}\alpha}B(\beta(w))\to k,
\end{equation*}
where these functions are homogeneous and their homogeneous components  are  zero by def\/inition unless
$a_w\in B_{\alpha(w)-1-}(\beta(w),A) $.

\subsubsection{Decorations in the string topology case}

The above procedure  is tantamount to changing the decorations as follows:

\vspace{-1mm}

\begin{table}[h]\centering \caption{Weights for open/closed string topology.}
\vspace{1mm}

\begin{tabular}{lll}
Type&Decoration& Weight $w(s)$\\
\hline
side without marked point\tsep{1pt}\\
\quad  part of ``in'' boundary&$a\in A_{\varnothing}$&$a$\\
\quad part of ``out''  boundary and part of a rectangle&$a\in A_{\varnothing}$&$a$\\
 \quad part of ``out''  boundary not part of a rectangle&$1\in A_{\varnothing}$&$1$\\
side with marked point marked by $\varnothing$&$a\in A_{\varnothing}$&$a$\\
side with marked point marked by $S$&$(a_S^1,a_S^2)$, $a_S^i\in A_{S}$&
$r^{\dagger}_{\varnothing S}(a_S^1 a_S^2)$
\end{tabular}\vspace{-1mm}
\end{table}

In \cite{hoch2} we used angle markings to this ef\/fect. In that language
the table above is the analog of the decorations for the
angle markings in the PROP case. The angles correspond exactly to the
components of the boundary minus the arcs.

\subsection{The action}

\begin{thm}
\label{mainthm} The open/closed $\beta$ brane labelled open/closed
 dg-PROP cell model of $\wSull$ provided by $\Sull_1$ acts
in a brane labelled open/closed  dg-PROP fashion on the brane
 labelled Hochschild complexes for a $\beta$-Frobenius algebra which
 satisfies
the Euler condition $(E)$.
\end{thm}

This has an immediate geometric consequence by using $\beta$ Frobenius
algebras coming form the geometric data of Section~\ref{geodatasection}.
Notice that in this case the bar complex
$B(\varnothing,\varnothing))$ after dualizing computes  $H_*(LM)$ where $LM$ is
the free loop space \cite{J} and the bar complex  $B(b,b')$ after dualizing
computes $H_*(PM,N_b,N_{b'})$ where $PM(N_b,N_{b'})$ is the space
of paths which start in $N_b$ and end in $N_{b'}$.

\begin{cor}
If $M$ is a simply connected compact manifold with a given set of
 $D$-branes realized by submanifolds $N_b:b\in\calB$ then there are open
 closed
string topology type operations on $H_*(LM)$ and the various $H_*(PM,N_b,N_{b'})$.
\end{cor}

Here def\/ine new operations through the $E_2$ term of the respective spectral sequence.
One actually can do it for the $E_1$ term; confer \cite{hoch2}.
This is what we call string topology type operations. One can ponder if such
operations exist in a purely geomtric framework and if these coincide with the ones def\/ined here.

\begin{figure}[t]
\centerline{\includegraphics[scale=0.96]{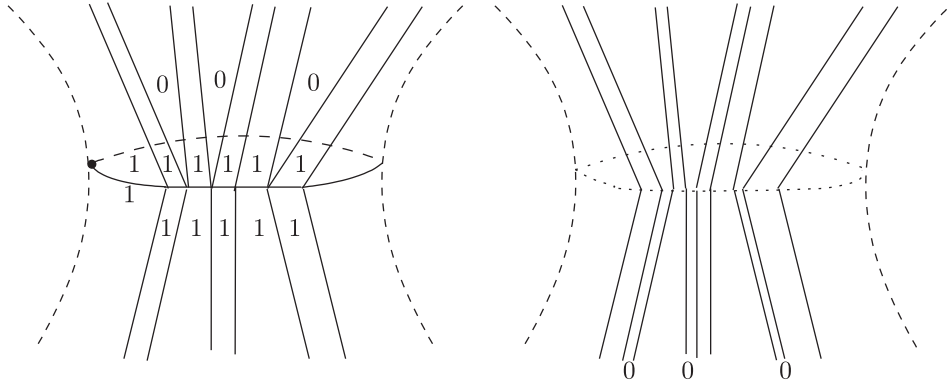}}
\caption{Gluing with angle labels in \cite{hoch2}, in the current
terminology  the label 0 corresponds
to a side without marked point on an ``out'' boundary that is not part of a rectangle.}\label{gluingcompat}
\vspace{-2mm}
\end{figure}

\begin{proof}[Proof of Theorem \ref{mainthm}]
There are two things left to prove, the dg-properties and the
 compati\-bi\-lity
of the gluings on both sides of the action.
We f\/ist show the dg-properties.
In the closed case the argument is the same as in \cite[Section 4.2.1]{hoch2}.
It basically relies on the fact that the co-multiplication and
the multiplication are dual to each other.  The argument carries over
 to the open case upon noticing that the only relevant case is the one where we locally
consider an arc which has no parallel arc, as otherwise the terms of the dif\/ferential cancel out.
The only new case is when this arc is the only arc in the window.
 In this case (unlike in the closed case) removing the incident arc will still leave
 two decorations. In all other cases one of the decorations before removing the arc will be by $1$.
In this new case there are however two decorations before
and after removing the arc due to the new rules of decorating on
open windows, and hence the right hand side of the equations~(4.7) and~(4.8) of
\cite{hoch2} still give the correlators on the surface
with the arc removed.

The argument that the gluings on the Hochschild side in view
of Theorem~\ref{discretethm} and the geometric side coincide is completely analogous to
\cite[Theorem 4.4]{hoch2}.
There are two steps in the argument. The f\/irst is that on the cellular chain side, only
the graphs with maximal dimension appear and we have to check that only these
appear on the Hochschild side. This is forced by the labelling. This labelling is equivalent
to considering the gluing of~\cite{hoch2} for angle labelled graphs (see Fig.~\ref{gluingcompat}).
If does not want to take this detour one can prove this directly by noticing that the degeneracies
duplicate edges on gluing. Now with this gluing the number of arcs is always maximal after gluing.
On the Hochschild side, this is automatic as the number
of ``in'' variables has to be equal to the number of ``out'' variables in order to obtain a non-zero
gluing. The second step is to check
that the actions given by the discretized arcs coincide before and after gluing. This is the content of
Theorem~\ref{discretethm}. The surprising new feature is that among the new local gluings coming
from the open sector the compatibility holds only if additionally the condition $(E)$ is satisf\/ied.
\end{proof}

\subsubsection{Examples multiplication and comultiplication\\ in the  open sector in the string topology case}
\label{sullexpar}

As an example we will consider the arc graphs given in Fig.~\ref{openmult} but now with the Sullivan decoration.
Given
\[
a=a^{(2)}_T\otimes a_1\oto a_n\otimes a^{(1)}_S\in B_n(A_T,A,A_S)
\]
and
\[
a'=a_S^{(2)}\otimes a'_1\oto a'_m\otimes a^{(1)}_U\in
B_m(A_S,A,A_U)
\]
 their product now lies
in $B_{n+m}(A_T,A,A_U)$ and is given by
\begin{equation*}
m_{TSU}(aa') = \bigg(\int_Sa^{(1)}_Sa^{(2)}_S\bigg)  a_T^{(2)}\otimes a_1\oto a_n \otimes a'_1\oto a'_m\otimes a^{(1)}_U.
\end{equation*}
The calculation goes as follows. The integrals are $\int
r^{\dagger}_T (a^{(1)}_T a^{(2)}_T)=\int_T a^{(1)}_T a^{(2)}_T$
which dualizes to ${\rm id}_T(a_T^{(1)})=a_T^{(1)}$ and likewise for $U$.
For the rectangles we get either $\int a_i a''_j$ or $\int a_i'
a''_k$ which dualize to ${\rm id}(a_i)$ and ${\rm id}(a_j)$. Finally we get
$\int r_S^{\dagger}(a_S^{(1)}a_S^{(2)})1=\int_S a^{(1)}_Sa^{(2)}_S$
which gives the factor in front.

The corresponding  co-products  $\Delta_{TSU}: B(A_T,A,A_U)\to B(A_U,A,A_S) \otimes B(A_S,A,A_T)$ are left unchanged.

\section{Moduli space actions}
\label{modulisection}

\begin{figure}[t]
\centerline{\includegraphics{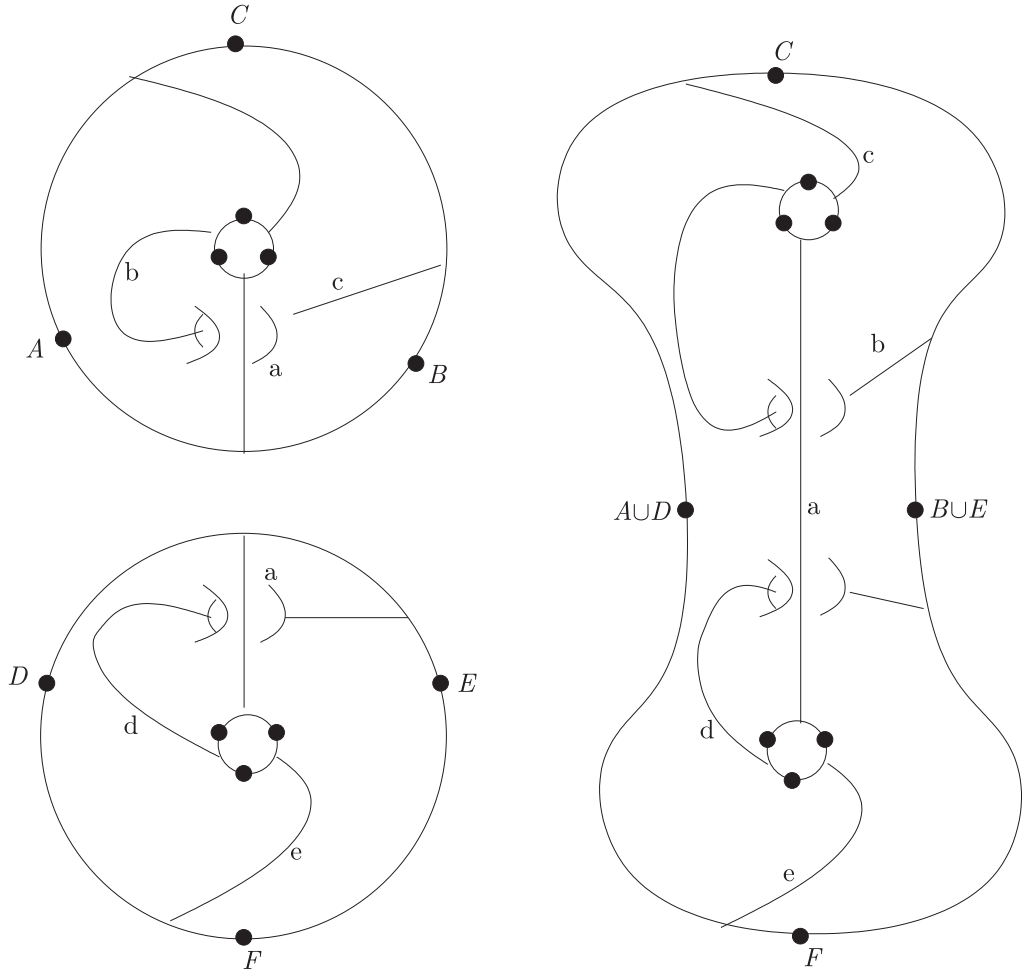}}
\caption{A gluing of two elements of moduli space
whose result does not lie in the moduli space. The window $AB$ is
glued to window $DE$ and the resulting element is not quasi-f\/illing
anymore.}\label{example}
\end{figure}

\begin{figure}[t]
\centerline{\includegraphics[scale=0.96]{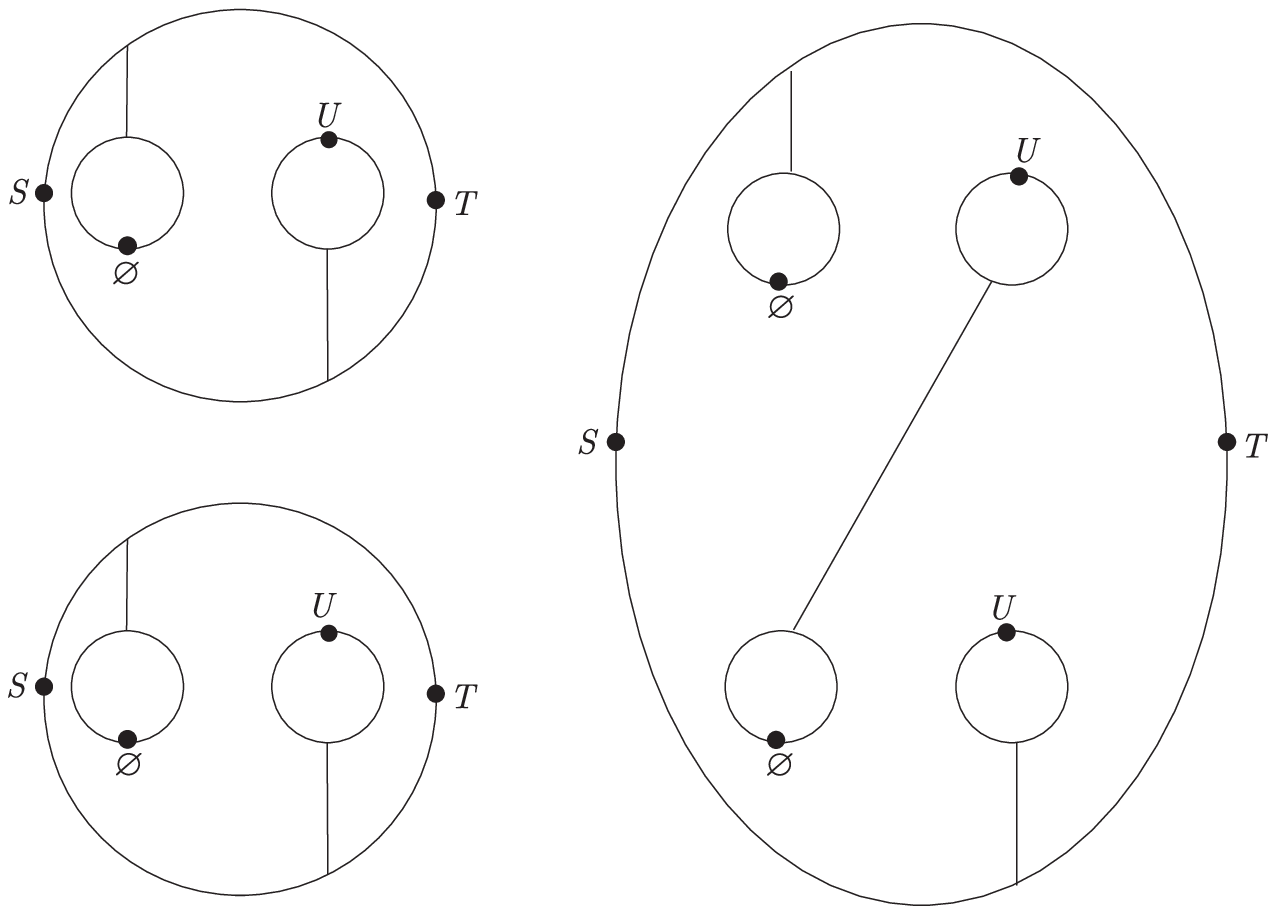}}
\caption{Gluing of two cells of moduli space whose
result does not lie in the moduli space. The gluing is along the
windows $ST$ which face each other. Again the result is not
quasi-f\/illing.}\label{example2}
\vspace{-2mm}
\end{figure}

In the moduli spaces case, there an associated chain complex indexed
by graphs. The dif\/f\/iculties in this case are manifold. First the
cells of the chain complex are open cells. As we saw in~\cite{hoch1,hoch2}, the way to deal with this is to pass to the
associated graded complex and look at the actions induced from the
topological level there. Now there are new problems that arise in
the open/closed case. While in the purely closed situation, the
``forbidden'' gluings which on the topological level gave rise to
elements outside the moduli space were codimension one, here there
are ``forbidden'' gluings for whole cells. An example is given in
Figs.~\ref{example} and~\ref{example2}. We can deal with this by
restricting the cells to come from an open/closed duality subspaces
of moduli space. Also there are problems since we have to deal with
internal marked points.

\subsection{Point clusters}
The property of moduli space that only once punctured polygons may
appear among the complementary regions  is very fragile under
gluing. Not even the inclusion of the open and closed sectors into
each other is stable with respect to this condition. So we will
allow polygons with an arbitrary number of punctures. Using Strebel
dif\/ferentials, we can put a conformal structure on such a piece and
have one distinguished point for this polygon. We now choose the
geometric interpretation that all the points in the polygon with
their brane labels are a cluster of points located at that
distinguished point. We can think of these points as labelled
``bosons'' sitting on top of each other. In this sense they are
equivalent to one internal marked point with a multiple brane label.
This is very close to the moduli spaces considered in~\cite{LM}.

\subsection{Open/closed duality}
\begin{figure}[t]
\centerline{\includegraphics[scale=0.96]{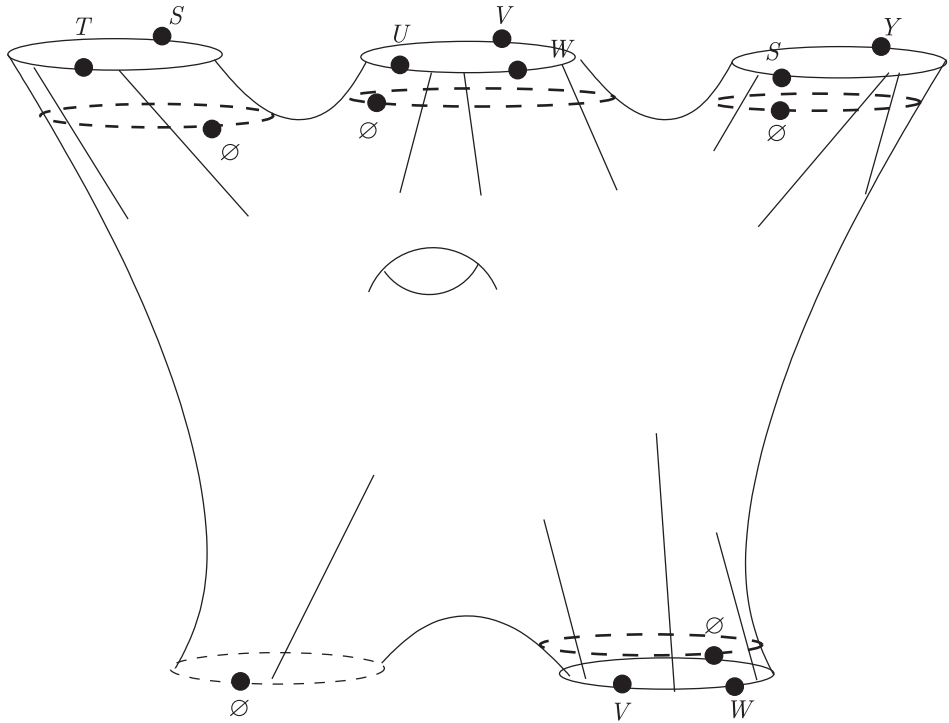}}
\caption{Open/closed duality: cut on the dotted lines.}\label{openclosed}
\vspace{-2mm}
\end{figure}

In \cite{KP} we saw that the open/closed duality holds on the chain level. This basically means
that any element of $\widetilde{\rm Arc}$ can be decomposed into a piece which is purely closed
and cylinders whose one end is closed. This is simply obtained by cutting on a closed
curve parallel to each boundary. If the boundary is already closed there is no need to cut.
An example of this is given in Fig.~\ref{openclosed}.
Notice the since we are cutting, we have the free choice of a point on the boundary and
hence the pieces are not unique. They are of course unique up to twisting on the boundary,
with the one parameter family on the annulus with closed windows, which is equivalent to moving a marked point marked by~$\varnothing$ along the boundary.
On the chain level this does however not hold for the whole space and hence we will impose this condition.

\begin{df}
We say that two arc families $\a$ and $\a'$
 are in {\em general position} with respect to the windows $w$ of $\a$ and $w'$ of $\a'$
 if the number of  arcs after the gluing is subadditive.
 \end{df}
 A main result of \cite{hoch1} is that elements of open cells of  two purely closed quasi-f\/illing arc graphs
 are not in general position only in codimension one. This was enough to induce an operad
 structure on the associated graded of the chains.

  Now, we also have  problems with gluings of the type where two f\/lags
 of a window belong to the same complementary region, e.g.\ the case~e)~2) discussed in Section~\ref{twocasespar}.

\begin{df}
 We say an open window of an arc family $\a$ is {\em degenerate}
 if its f\/lags are both part of the boundary  of a complementary region, but are not
 an edge of this boundary. The latter can only happen if there
 is only one marked point on the boundary. This case will be non-degenerate.

An {\em arc family is non-degenerate} if it has no degenerate
windows.
\end{df}

\subsubsection{Cell complexes}
For a $\beta$ arc graph $\alpha$ on $F$, we let $C(\alpha)$ be the set of
all projective weightings on $\alpha$. Recall that a weighting is by positive reals,
so that this set is just the open inside of a cell.

The space ${\rm Arc}(F,\beta)(n,m):=\widetilde{{\rm Arc}}(F,\beta)(n,m) /\R_{>0}$ has decomposition into these  open cells
\begin{equation*}
{\rm Arc}(F,\beta):=\sqcup_{\Gamma : \text{$\beta$ arc graph on $F$ with only active windows}} \dot C(\alpha).
\end{equation*}
This gives a complex whose generators are the open cells and whose boundary
is given by the dif\/ferential in the associated CW complex.
The terms in the sum are only over those graphs which appear in the boundary. These
are the subgraphs with one fewer edge.

\begin{df}
The open/closed duality moduli space is the space
given by those arc families which satisfy the following conditions
\begin{enumerate}\itemsep=0pt
\item The complementary regions are  only polygons possibly with punctures of any f\/inite number. \item The arc family admits a
  decomposition under the open/closed duality as above such that
  \begin{enumerate}\itemsep=0pt
  \item the pieces are in general position with respect to the boundaries obtained by cutting and
  \item  the annuli appearing in the decomposition are non-degenerate.
\end{enumerate}
\end{enumerate}

In particular we let  $\Mco$ be the component where the arc families are on $(F_{g,n}^{s},\beta)$
where the point clusters are given by the sets of marked points of cardinality $s_i$ and $s=\sum_{s_i}$, where a point
cluster is the set of marked points within one polygonal complementary region.
\end{df}

Notice that the spaces $\Mco$ are stratif\/ied by the spaces $\Mcos$ and each space
 $\Mcos$ is a f\/inite unramif\/ied cover of $\M_{g,\delta_1,\dots,\delta_n}^k$.

To avoid yet additional notation, we will think of $\M_{g,\delta_1,\dots,\delta_n}^s$ and $\Mco$ as a subspace of ${\rm Arc}(F,\beta)(n_1,n_2)$ where $n_1$ is the number of $\delta_i=1$ and $n_1+n_2=n$ and $F=F_{g,n}^s$.

The subspaces are then just given by the disjoint union of open cells of those graphs that
satisfy the additional requirements.

There are gluings on the topological level, which are giving by scalings. Given $\a$ and a~window~$w$ of it together with $\a'$ and a window $w'$ on it, we scale all arcs $\a$ by $\a'(w')$ and all arcs of~$\a'$ by~$\a(w)$, just as in~\cite{KLP}.  After the scaling the two windows have the same weight and
we can glue. The structure we get is a two colored operad (open/closed) with the additional
information of brane labels. We call such a structure a brane labelled open/closed operad.

\begin{lem}
This gluing yields a brane-labelled open/closed operad structure. And this induces
an operad structure on the complex of open cells.
\end{lem}

\begin{proof}
On the topological level the only thing that is left to be checked is the associativity. Adapting \cite{KLP} this is straightforward.
For the chains, we notice that the set obtained from composing cells is a union of cells and proceed
as in \cite{hoch1}.
\end{proof}

If we only stick to basic brane labels, we would need to introduce more colors, which would be
as usual pairs $(S,T)$ of $S,T\neq \varnothing$ or $(\varnothing,\varnothing)$.

\begin{lem}
The open cell complexes of ${\rm Arc}(F,\beta)(n_1,n_2)$  are graded by dimension~-- which for~$\dot C(\Gamma)$ is the number of arcs minus one~-- and
the induced operad structure respects the corresponding filtration and hence passes to the associated
graded complexes.
\end{lem}
\begin{proof}
This follows from the fact that when gluing two windows with $k$ and $l$ arcs, the maximal
number is $k+l+1$.
\end{proof}

\begin{prop}
The associated graded of the complexes of open cells  $\Mco$ form a~suboperad of
the associated graded of the complexes of open cells of ${\rm Arc}(F,\beta)(n_1,n_2)$.
\end{prop}

\begin{proof}
We have to show that when gluing two  open cells which are in the open/closed moduli part,
we are still in the open/closed moduli part up to codimension one. For the closed part this is contained in~\cite{hoch1} up to handling the punctures. But since we allow arbitrary punctures in polygons, we do not have to keep track of them. This generalizes directly to gluing on two closed windows in the possible presence of additional open windows. Now we use the condition of the open/closed moduli space to decompose
any given two elements into cylinders and pieces that only have closed boundaries.
We now only have to consider gluing on such a cylinder along an open window to an open window.
Since the cylinder will have a non-degenerate quasi-f\/illing arc family, we see that indeed the
result of the gluing will be quasi-f\/illing and non-degenerate as we will locally never glue
polygons twice on sides that are not consecutive.
\end{proof}

To each such cell $\dot C(\Gamma)$ we can associate the correlation function $Y(\Gamma)$ acting
on the approp\-riate Hochschild co-chains or bar complexes. Now
since we passed to the associated graded on the moduli space side, we will have to do the
same thing on the algebraic side. This is accomplished by grading with respect to the number of comultiplications, analogous to~\cite{hoch2} and then passing to the associated graded.
The resulting objects can naturally be called a $\beta$-labelled $\rm Hom$ operad.

\begin{thm}
\label{mainthmtwo}
There is  an operadic cell model
associated to the $\beta$-brane labelled open/closed
moduli spaces $\Mco$ which acts on $\beta$-labelled Hochschild co-chains via operadic correlation functions with values in a $\beta$-labelled $Hom$ operad.
\end{thm}

\begin{rmk}
On the image this operation is dg with respect to the induced dif\/ferential.
\end{rmk}

\section{Outlook}
Again like in \cite{romp,postnikov} we can consider the stabilization. We see that in this
case, we need that all the Frobenius algebras $A_S$ are (normalized) semi-simple in
order to pass to the appropriate stabilization.

One case where this would be true would be in Landau--Ginzburg models.
We are currently working on the details of this theory.

One can furthermore ask about the modular operad structure on the moduli space.
Then further technical complications arise from the intricate structure of the f\/lows def\/ining
the chain level structure of the c/o structure in~\cite{KP}.  In this case we will show
that there is an underlying solution to the quantum master equation.

\appendix

\section{Appendix: operadic, PROPic and c/o structures}\label{appendixA}
% use *-form to suppress numbering

\subsection{The def\/inition of a c/o structure}
\label{coapp}
 Specify an object ${\mathcal O}(S,T)$ in some f\/ixed
symmetric monoidal category for each pair $S$ and $T$ of f\/inite
sets.  A {\it $G$-coloring} on ${\mathcal O}(S,T)$ is the further
specif\/ication of an object $G$ in this category and a morphism
$\mu : S\sqcup T \to {\rm Hom}({\mathcal O}(S,T),$G$)$, and we shall let
${\mathcal O}_\mu (S,T)$ denote this pair of data.

A $G$-colored ``closed/open'' or {\it c/o structure} is a
collection of such objects ${\mathcal O}(S,T)$ for each pair of
f\/inite sets $S$, $T$ together with a choice of weighting $\mu$ for
each object supporting the following four operations which are morphisms in the
category:

\vskip .1in

\noindent
\noindent{\it Closed gluing}:  $\forall\, s\in S$, $\forall\, s'\in S'$ with $\mu (s)=\mu '(s')$,
\[
\circ_{s,s'}: \  {\mathcal O}_\mu (S,T)\otimes{\mathcal O}_{\mu '} (S',T')\to {\mathcal O}_{ \mu ''}
(S\sqcup S'-\{ s,s'\}
, T\sqcup T');
\]

\noindent{\it Closed self-gluing}: $\forall \, s,s'\in S$ with $\mu (s)=\mu (s')$ and $s\neq s'$,
\[
\circ^{s,s'}: \ {\mathcal O}_\mu (S,T)\to {\mathcal O}_{\mu ''} (S-\{ s,s'\} , T);
\]

\noindent{\it Open gluing}:
$\forall  \, t\in T$, $\forall \, t'\in T'$ with $\mu (t)=\mu '(t')$,
\[
\bullet_{t,t'}: \ {\mathcal O}_\mu (S,T)\otimes{\mathcal O}_{\mu '} (S',T')\to {\mathcal O}_{\mu ''}
(S\sqcup S' , T\sqcup T'-\{t,t'\});
\]

\noindent{\it Open self-gluing}: $\forall \, t,t'\in T$ with $\mu (t)=\mu (t')$ and $t\neq t'$,
\[
\bullet^{t,t'}:  \ {\mathcal O}_\mu (S,T)\to {\mathcal O}_{\mu ''} (S , T-\{ t,t'\}).
\]

\noindent In each case, the coloring $\mu ''$ is  induced in the
target in the natural way by restriction, and we assume that $S\sqcup S'\sqcup T\sqcup T'-\{
s,s',t,t'\}\neq\varnothing$.

The  axioms are that the operations are equivariant
for bijections of sets and for bijections of pairs of sets, and the collection of all operations taken together
satisfy associativity.

Notice that we use the formalism of operads indexed by f\/inite
sets rather than by natural numbers as in \cite{MSS} for instance.

\subsection{Restrictions}
A c/o structure specializes to standard algebraic objects
in the following several ways.

There are the two restrictions
 $(\CO(S,\varnothing),\circ_{s,s'})$ and
 $(\CO(\varnothing,T),\bullet_{\t,\t'})$
each of which forms a~$G$-colored cyclic operad in the usual sense.

The spaces
$(\CO(S,T), \circ_{s,s'},\bullet_{\t,\t'})$ with only the non self-gluings as structure maps
form a~cyclic
$G\times \Z/2\Z$-colored operad, where the $\Z/2\Z$ accounts for
open and closed, e.g., the windows labeled by $S$ are
regarded as colored by $0$ and the windows labeled by $T$ are
regarded as colored by~$1$.

If the underlying category has a
coproduct (e.g., disjoint union for sets and topological
spaces, direct sum for Abelian groups and linear spaces), which we
denote by $\coprod$, then the indexing sets can be regarded as
providing a grading: i.e.,   $(\coprod_T
\CO(S,T),\circ_{s,s'})$ form a cyclic $G$-colored operad graded by
the sets $T$, and
$(\coprod_S\CO(S,T),\bullet_{\t,\t'})$ form a cyclic
$G$-colored operad graded by the sets $S$.

\subsection{Modular properties}
\label{modular}
There is a relationship between c/o structures and modular operads. Recall
that in a modular operad there is an additional grading on the
objects, which is additive for gluing and increases by one for
self-gluing. Imposing this type of grading here, we
def\/ine a {\it $(g,\chi -1)$ c/o structure} to be a c/o structure with
two gradings $(g,\chi)$,
\[
\CO(S,T)=\coprod_{g\geq 0,\chi
\leq 0} \CO(S,T;g,\chi)
\] such that
\begin{enumerate}
\item[(1)]
$\CO(S,T;\chi-1)=\coprod_{g\geq 0}\CO(S,T;g,\chi)$ is additive in
$\chi-1$ for $\bullet_{t,t'}$, and $\chi-1$ increases by one for
$\bullet^{t,t'}$; and

\item[(2)] $\CO(S,T;g)=\coprod_{\chi\leq
0}\CO(S,T;g,\chi)$ is additive in $g$ for $\circ_{s,s'}$, and $g$
increases by one for $\circ^{s,s'}$.
\end{enumerate}

It follows that a $(g,\chi -1)$
c/o structure is a modular $G$-colored bi-operad in the sense that the
 $\CO(S,T;g)$ form a $T$-graded $\Rp$-colored
 modular operad\footnote{We impose
neither  $3g-3+|S|>0$ nor $3(-\chi+1)+|T|-3\geq 0$.} for the  gluings $\circ_{s,s'}$ and
$\circ^{s,s'}$, and the $\CO(S,T;1-\chi)$ form an $S$-graded
$\Rp$-colored  modular operad$^*$ for the gluings $\bullet_{t,t'}$ and
$\bullet^{t,t'}$.

\subsection{Brane-labeled c/o structures}

A  brane-labeled c/o
structure is a c/o structure $\{\CO(S,T)\}$  together with a f\/ixed
Abelian monoid~$\mathcal{P}$ of brane labels and for each $\a\in
\CO(S,T)$ a bijection $N_{\a}:T\rightarrow T$ and a bijection
$(\lambda_{\a},\rho_{\a}): T\rightarrow \mathcal{P} \times
\mathcal{P}$, such that
\begin{itemize}\itemsep=0pt
\item[(1)]
 $\rho(t)=\lambda(N(t))$,
 \item[(2)] if
$N_{\a}(t)\neq t$ and $N_{\a'}(t')\neq t'$
\begin{alignat*}{3}
& N_{\a\bullet_{t,t'}\a'}(N_{\a}^{-1}(t))=N_{\a'}(t'), \qquad &&
N_{\a\bullet_{t,t'}\a'}(N_{\a'}^{-1}(t'))=N_{\a}(t), & \\
& \rho_{\a\bullet_{t,t'}\a'}(N^{-1}(t))=\lambda_{\a}(t)\rho_{\a'}(t'),\qquad &&
\lambda_{\a\bullet_{t,t'}\a'}(N(t))=\lambda_{\a'}(t')\rho_{\a}(t), &
\end{alignat*}

\item[(3)] $N_{\a}(t)\neq t$ and $N_{\a'}(t')\neq t'$
\begin{alignat*}{3}
& N_{\bullet^{t,t'}(\a)}(N_{\a}^{-1}(t))=N_{\a}(t'),\qquad &&
N_{\bullet^{t,t'}(\a)}(N_{\a}^{-1}(t'))=N_{\a}(t), & \\
& \rho_{\a\bullet^{t,t'}\a'}(N^{-1}(t))=\lambda_{\a}(t)\rho_{\a}(t'),\qquad
&&\lambda_{\a\bullet^{t,t'}\a'}(N(t))=\lambda_{\a}(t')\rho_{\a}(t),&
\end{alignat*}

\item[(4)] if either $N_{\a}(t)=t$ or $N_{\a'}(t')=t'$ but not
both, then in the above formulas, one should substitute
 $N_{\a'}(t')$ for $N_{\a}(t)$ in the f\/irst case and inversely in
 the second case. (If both $N_{\a}(t)=t$ and
 $N_{\a'}(t')=t'$, then there is no equation.)
\end{itemize}

This is the axiomatization  of the geometry given by open windows
with endpoints labeled by right ($\rho$) and left ($\lambda $) brane labels, their order and
orientation along the boundary components induced by the
orientation of the surface, and the behaviour of this data under
gluing.

 For a brane-labeled c/o structure and an idempotent submonoid $\B
 \subset \mathcal{P}$ (i.e., for all $b\in \B$, $b^2=b$),
 one has the {\rm
$\B\times \B$-colored} substructures def\/ined by restricting the
gluings~$\bullet_{t,t'}$ and~$\bullet^{t,t'}$ to compatible colors
$\lambda(t)=\rho(t')$.

\subsection{PROPs and partial modular operads}

Recall that a PROP is a collection of objects $\calO(S,T)$ as above
with two sets of operations. Composition: for each bijection
$T\leftrightarrow S'$ a morphism $\circ_{\phi}:\calO(S,T)\otimes
\calO(S',T')\to \calO(S,T')$. We can think of $S$ as ``in'' labels
and $T$ as ``out'' labels. Mergers or vertical compositions:
$\calO(S,T)\otimes \calO(S',T')\to \calO(S\amalg S',T\amalg T)$
satisfying the obvious compatibilities, associativities and
invariance under bijections.

A partial modular operad is a modular operad in which not all gluing
operations need to be def\/ined, when they are def\/ined all the
relations hold. A particular kind of partial modular  operad is a
colored modular operad. There is an additional coloring of $S$ by
some set $C$ and only like colors can be glued. Another partial
structure is the following: For each element in~$\calO(S)$ f\/ix   a
partition if $S$ into $S_0\amalg S_1$. We can think of $S_0$ as
``in'' and $S_1$ as ``out''. The gluings are given for
$\circ_{s,s'}:\calO(S)\otimes \calO(S')\to
\calO(S\setminus\{s\}\amalg S'\setminus\{s'\})$ for all $s\in
S_1$, $s'\in S'_0$. We will call this an i/o modular operad.

Notice and i/o modular operad in a category with a direct sum induces a PROP.

Given a bijection $\phi: T\to S'$ we  perform the concatenation of
all operations $\circ_{t,\phi(t)}$, where the f\/irst gluing is a
non-self gluing and all further gluings are self-gluings. The
mergers are then just def\/ined as the coproduct, which is the
disjoint union.

\subsection{Two colored or open/closed colored PROPs}

A two colored or open/closed PROP is a PROP with a partition for
each the sets $S$, $T$ that is it is given by objects $\calO(S_o\amalg
S_c,T_o\amalg T_c)$ and there are gluing operations for bijections
$S_o \leftrightarrow T'_o$ and $S_c \leftrightarrow T'_c$ that give
PROP structures separately (treating the other set as a grading) and
jointly.

\subsection{Weaker structures on the topological level}
\label{quasiproppar}

A topological quasi-PROP has the data of PROP
but the associativity only needs to hold up to homotopy. Notice that this is
enough to guarantee that there is a PROP structure on the homology
level. The adjective ``quasi'' in any context means that the
associativity need not hold and the specif\/ication ``topological
quasi'' means that the objects are  topological spaces and the
associativity holds up to homotopy.

\subsubsection{C/o versions, c/o PROP}
\label{coprop}
To get the c/o versions, we need two additional colorings on top of all the
other structures. This first is open/closed and the second is a coloring
$\mu$ by $\mathbb R$. The gluing then is partial or colored with respect to
both the coloring open/closed and the coloring $\mu$.

Thus a c/o PROP has colors ``closed'' and ``open'' together with a
real weight on each ``in'' and each ``out''. There will be two types
of full gluings: closed and open given by bijections between the
open ``ins'' and ``outs'' or the closed ``ins'' and ``outs''. The
gluings are possible if the bijection $\phi$ has the property that
the weight of $t$ is the weight of $\phi(t)$. Both gluings and the
simultaneous open/closed gluing will satisfy the usual equations.

Adding brane labels is straightforward as above.

\section*{Acknowledgments}
This project was completed in many phases over the last couple of years.
We wish to thank the institutes at which we worked on it for their hospitality and support. These are the
MPI, IHES where this work was conceived, and the CTQM in Aarhus where many of the f\/inal
details were f\/ixed. We also wish to thank Bob Penner and J{\o}rgen Andersen for discussions.
Special thanks goes to the organizers of the conference ISQS 2009 which
provided the forum to present these ideas and prompted me to f\/inish
the project and write up the results.

We also gratefully acknowledge support from the NSF grant DMS-080588.

\pdfbookmark[1]{References}{ref}

\LastPageEnding


\begin{thebibliography}{99}

\footnotesize\itemsep=0pt

\bibitem{ramirez}
Baas N.A., Cohen  R.L., Ram\'{\i}rez A.,
The topology of the category of open and closed strings,
in Recent Developments in Algebraic Topology,
{\it Contemp. Math.}, Vol.~407, Amer. Math. Soc., Providence, RI, 2006,  11--26,
\href{http://arxiv.org/abs/math.AT/0411080}{math.AT/0411080}.

\bibitem{BCT}
Blumberg A.J., Cohen R.L., Teleman C.,
Open-closed f\/ield theories, string topology, and Hochschild homo\-lo\-gy,
in Alpine Perspectives on Algebraic Topology, Editors C.~Ausoni, K.~Hess and J.~Scherer,
{\it Contemp. Math.}, Vol.~504, Amer. Math. Soc., Providence, RI, 2009, 53--76,
\href{http://arxiv.org/abs/0906.5198}{arXiv:0906.5198}.

%\bibitem{CJ}
%Cohen R.L., Jones J.D.S.,
%A homotopy theoretic realization of string topology,
%\href{http://dx.doi.org/10.1007/s00208-002-0362-0}{{\it Math. Ann.}} {\bf 324} (2002), 773--798,
%\href{http://arxiv.org/abs/math.GT/0107187}{math.GT/0107187}. !!!Not cited!!!

 \bibitem{CS}
 Chas M., Sullivan D.,
String topology,
{\it Ann. of Math.}, to appear,
\href{http://arxiv.org/abs/math.GT/9911159}{math.GT/9911159}.

\bibitem{godin}
Godin V.,
Higher string topology operations,
\href{http://arxiv.org/abs/0711.4859}{arXiv:0711.4859}.

\bibitem{HVZ}
Harrelson E., Voronov A.A., Zuniga J.J.,
Open-closed moduli spaces and related algebraic structures,
\href{http://arxiv.org/abs/0709.3874}{arXiv:0709.3874}.

\bibitem{J}
Jones J.D.S.,
Cyclic homology and equivariant homology,
\href{http://dx.doi.org/10.1007/BF01389424}{{\it Invent. Math.}} {\bf 87} (1987), 403--423.


%\bibitem{KLi1}
%Kapustin A., Li Y.,
%D-branes in Landau--Ginzburg models and algebraic geometry,
%\href{http://dx.doi.org/10.1088/1126-6708/2003/12/005}{{\it J. High Energy Phys.}} {\bf 2003} (2003), no.~12, 005, 44~pages,
%\href{http://arxiv.org/abs/hep-th/0210296}{hep-th/0210296}. !!!Not cited!!!

%\bibitem{KR}
%Kapustin A., Rozansky L.,
%On the relation between open and closed topological strings,
%\href{http://dx.doi.org/10.1007/s00220-004-1227-z}{{\it  Comm. Math. Phys.}} {\bf  252} (2004), 393--414,
%\href{http://arxiv.org/abs/hep-th/0405232}{hep-th/0405232}. !!!Not cited!!!

\bibitem{hoch1}
Kaufmann R.M.,
Moduli space actions on the Hochschild co-chains of a Frobenius algebra. I.~Cell operads,
\href{http://dx.doi.org/10.4171/JNCG/10}{{\it J. Noncommut. Geom.}} {\bf 1} (2007), 333--384,
\href{http://arxiv.org/abs/math.AT/0606064}{math.AT/0606064}.

\bibitem{hoch2}
Kaufmann R.M.,
Moduli space actions on the Hochschild co-chains of a Frobenius algebra. II.~Correlators,
\href{http://dx.doi.org/10.4171/JNCG/22}{{\it J. Noncommut. Geom.}} {\bf 2} (2008), 283--332,
\href{http://arxiv.org/abs/math.AT/0606065}{math.AT/0606065}.


\bibitem{romp}
Kaufmann R.M.,
Noncommutative aspects of open/closed strings via foliations,
\href{http://dx.doi.org/10.1016/S0034-4877(08)80016-X}{{\it Rep. Math. Phys.}} {\bf 61} (2008), 281--293,
\href{http://arxiv.org/abs/0804.0608}{arXiv:0804.0608}.

\bibitem{postnikov}
Kaufmann R.M.,
Dimension vs. genus: a surface realization of the little $k$-cubes and an $E_{\infty}$ operad,
in Algebraic Topology -- Old and New (M.M.~Postnikov Memorial Conference),
\href{http://dx.doi.org/10.4064/bc85-0-17}{{\it Banach Center Publ.}}, Vol.~85, Polish Acad. Sci., Warsaw, 2009, 241--274,
\href{http://arxiv.org/abs/0801.0532}{arXiv:0801.0532}.

\bibitem{KP}
Kaufmann R.M., Penner R.B.,
Closed/open string diagrammatics,
\href{http://dx.doi.org/10.1016/j.nuclphysb.2006.03.036}{{\it Nuclear Phys.~B}} {\bf 748} (2006), 335--379,
\href{http://arxiv.org/abs/math.GT/0603485}{math.GT/0603485}.

\bibitem{KLP}
Kaufmann R.M., Livernet  M., Penner R.B.,
Arc operads and arc algebras,
\href{http://dx.doi.org/10.2140/gt.2003.7.511}{{\it Geom. Topol.}} {\bf 7} (2003), 511--568,
\href{http://arxiv.org/abs/math.GT/0209132}{math.GT/0209132}.

\bibitem{LM}
Losev A., Manin Yu.,
New moduli spaces of pointed curves and pencils of f\/lat connections,
\href{http://dx.doi.org/10.1307/mmj/1030132728}{{\it Michigan Math.~J.}} {\bf 48}  (2000), 443--472,
\href{http://arxiv.org/abs/math.AG/0001003}{math.AG/0001003}.
%Dedicated to William Fulton on the occasion of his 60th birthday.

\bibitem{MSS}
 Markl M., Shnider S., Stashef\/f J.,
Operads in algebra, topology and physics,
{\it Mathematical Surveys and Monographs}, Vol.~96, Amer. Math. Soc., Providence, RI, 2002.


%\bibitem{Merk}
%Merkulov S.A.,
%De Rham model for string topology,
%\href{http://dx.doi.org/10.1155/S1073792804132662}{{\it  Int. Math. Res. Not.}} {\bf  2004} (2004),  no.~55, 2955--2981,
%\href{http://arxiv.org/abs/math.AT/0309038}{math.AT/0309038}. !!!Not Cited in the paper!!!


\bibitem{P1}
Penner R.C.,
Decorated Teichm\"uller theory of bordered surfaces,
\href{http://projecteuclid.org/getRecord?id=euclid.cag/1098468019}{{\it Comm. Anal. Geom.}} {\bf  12}  (2004),  793--820,
\href{http://arxiv.org/abs/math.GT/0210326}{math.GT/0210326}.

\bibitem{Qu}
Quillen D.,
Elementary proofs of some results of cobordism theory using Steenrod operations,
\href{http://dx.doi.org/10.1016/0001-8708(71)90041-7}{{\it Adv. Math.}} {\bf 7} (1971), 29--56.

\bibitem{St}
Strebel K.,
Quadratic dif\/ferentials,
{\it Ergebnisse der Mathematik und ihrer Grenzgebiete (3)}, Vol.~5,
Springer-Verlag, Berlin, 1984.

\bibitem{C1}
Sullivan D.,
Sigma models and string topology,
in  Graphs and Patterns in Mathematics and Theoretical Physics,
{\it Proc. Sympos. Pure Math.}, Vol.~73, Amer. Math. Soc., Providence, RI, 2005,  1--11.


\bibitem{C2}
Sullivan D.,
Open and closed string f\/ield theory interpreted in classical algebraic topology,
in  Topology, Geometry and Quantum Field Theory,  \href{http://dx.doi.org/10.1017/CBO9780511526398.014}{{\it London Math. Soc. Lecture Note Ser.}}, Vol.~308,
Cambridge Univ. Press, Cambridge, 2004,  344--357,
\href{http://arxiv.org/abs/math.QA/0302332}{math.QA/0302332}.


\bibitem{sullsurv}
Sullivan D.,
String topology: background and present state,
\href{http://arxiv.org/abs/0710.4141}{arXiv:0710.4141}.

\bibitem{TZ}
Tradler T., Zeinalian M.,
Algebraic string operations,
\href{http://dx.doi.org/10.1007/s10977-007-9005-2}{{\it $K$-Theory}} {\bf 38} (2007), 59--82,
\href{http://arxiv.org/abs/math.QA/0605770}{math.QA/0605770}.

\end{thebibliography}
\end{document}